%% file: MDLin2014-arXiv.tex
\documentclass[11pt]{article}
\usepackage{epsf}
\usepackage{srcltx}
\usepackage{color}
\usepackage{float}
\usepackage{url}
\usepackage{enumerate,graphicx,subfig}
\usepackage{amsmath}
\input{MyMacros}
\newcommand{\be}{\begin{eqnarray}}
\newcommand{\ee}[1]{\label{#1}\end{eqnarray}}
\newcommand{\nn}{\nonumber \\}
\newcommand{\ese}{\end{eqnarray*}}
\newcommand{\bse}{\begin{eqnarray*}}
\newcommand{\rf}[1]{(\ref{#1})}
\newcommand{\half}{ \mbox{\small$\frac{1}{2}$}}
\def\TV{{\hbox{\scriptsize\rm TV}}}
\def\bZ{{\mathbf{Z}}}
\def\oomega{{\overline{\omega}}}

\def\oE{{\overline{E}}}
\def\oPsi{{\overline{\Psi}}}
\def\vixf{{{\rm vi}$(X,F)$ }}

\begin{document}

\title{Mirror Prox Algorithm for Multi-Term Composite Minimization and Semi-Separable Problems}
\author{Niao He\thanks{Georgia Institute of Technology, Atlanta, Georgia 30332, USA (\{nhe6,nemirovs\}@isye.gatech.edu). Research of these authors was supported by NSF Grants
CMMI-1232623 and CCF-1415498.}  \and Anatoli Juditsky\thanks{LJK, Universit\'e Grenoble Alpes,, B.P. 53, 38041 Grenoble Cedex 9, France (anatoli.juditsky@imag.fr). Research of this author was supported by the CNRS-Mastodons project GARGANTUA,
and the LabEx PERSYVAL-Lab (ANR-11-LABX-0025).} \and Arkadi Nemirovski\footnotemark[1]}
\maketitle


\begin{abstract}
In the paper, we develop a composite version of Mirror Prox algorithm for solving convex-concave saddle point problems and monotone variational inequalities of special structure, allowing to cover saddle point/variational  analogies of what is usually called ``composite minimization'' (minimizing a sum of an easy-to-handle nonsmooth and a general-type smooth convex functions ``as if'' there were no nonsmooth component at all). We demonstrate that the composite Mirror Prox inherits the favourable (and unimprovable already in the large-scale bilinear saddle point case) $O(1/\epsilon)$ efficiency estimate of its prototype. We demonstrate that the proposed approach can be successfully applied to Lasso-type problems with several penalizing terms (e.g. acting together $\ell_1$ and nuclear norm regularization) and to problems of semi-separable structures considered in the alternating directions methods, implying in both cases methods with the $O(1/\epsilon)$ complexity bounds.
\vspace{0.5cm}
\paragraph{Keywords:} numerical algorithms for variational problems, composite optimization, minimization problems with multi-term penalty, proximal methods
\vspace{0.5cm}
\paragraph{MSC(AMS):} 65K10, 65K05, 90C06, 90C25, 90C47.
\end{abstract}


\section{Introduction}
\subsection{Motivation}
Our work is inspired by the recent trend of seeking efficient ways for solving problems with hybrid regularizations or mixed penalty functions in fields such as machine learning, image restoration, signal processing and many others. We are about to present two instructive examples (for motivations, see, e.g., \cite{Buades05,Chambol05,Candes11}).

\paragraph{Example 1. (Matrix completion)} {Our first motivating example is matrix completion problem, where we want to reconstruct the original matrix $y\in\bR^{n\times n}$, known to be both sparse and low-rank, given noisy observations of part of the entries. Specifically, our observation is
$b=P_\Omega y+\xi$, where  $\Omega$ is a given set of cells in an $n\times n$ matrix, $P_\Omega y$ is the restriction of $y\in\bR^{n\times n}$ onto $\Omega$, and $\xi$ is a random noise.  A natural way to recover $y$ from $b$ is to solve the optimization problem}\begin{equation}\label{example: completion}
\Opt=\min\limits_{y\in\bR^{n\times n}}\left\{\frac{1}{2}\|P_\Omega y-b\|_2^2+\lambda\|y\|_1+\mu\|y\|_\nuc \right\}
\end{equation}
where $\mu,\lambda>0$ are regularization parameters. {Here  $\|y\|_2= \sqrt{\Tr(y^Ty)}$  is the Frobenius norm,  $\|y\|_1=\sum_{i,j=1}^n|y_{ij}|$ is the $\ell_1$-norm, and $\|y\|_\nuc=\sum_{i=1}^n\sigma_i(y)$ ($\sigma_i(y)$ are the singular values of $y$) is the nuclear norm of a matrix $y\in\bR^{n\times n}$.}

\paragraph{Example 2. (Image recovery)} {Our second motivating example is image recovery problem, where we want to recover an image $y\in\bR^{n\times n}$ from its noisy observations $b=Ay+\xi$, where $Ay$ is a given affine mapping (e.g. the restriction operator $P_\Omega$ defined as above, or some blur operator), and $\xi$ is a random noise. Assume that the image can be decomposed as $y=y_{\rm L}+y_{\rm S}+y_{\rm sm}$ where $y_{\rm L}$ is of low rank, $y_{\rm sm}$ is the matrix of contamination by a ``smooth background signal'', and $y_{\rm S}$ is a sparse matrix of ``singular corruption.'' Under this assumption in order to recover $y$ from $b$ it is natural to solve the optimization problem}
\begin{equation}\label{example: recovery}
\Opt=\min\limits_{y_{\rm L},y_{\rm S},y_{\rm sm}\in\bR^{n\times n}}\left\{\|A(y_{\rm L}+y_{\rm S}+y_{\rm sm}) -b\|_2+\mu_1\|y_{\rm L}\|_\nuc +\mu_2 \|y_{\rm S}\|_1+\mu_3 \|y_{\rm sm}\|_\TV\right\}
\end{equation}
where $\mu_1,\mu_2,\mu_3>0$ are regularization parameters. Here {$\|y\|_{\rm TV}$ is} the total variation of an image $y$:
$$
\begin{array}{c}
\|y\|_{\TV}=\|\nabla_i y\|_1+\|\nabla_j y\|_1,\\
\begin{array}{rcl}(\nabla_iy)_{ij}&=&y_{i+1,j}-y_{i,j},\; [i;j]\in\bZ^2:\;1\leq i<n-1,1\leq j<n,\\
(\nabla_jy)_{ij}&=&y_{i,j+1}-y_{i,j},\;[i;j]\in\bZ^2:\;1\leq i<n,1\leq j<n-1.
\end{array}
\end{array}
$$

{These and other examples motivate addressing} the following \emph{multi-term composite minimization problem}
\begin{equation}\label{model:unconstrained:0}
\min\limits_{y\in Y}
\left\{ \sum_{k=1}^K\left[\psi_k(A_ky+b_k)+\Psi_k(A_ky+b_k)\right]
\right\},
\end{equation}
and, more generally, the \emph{semi-separable problem}
\begin{equation}\label{model:constrained:0}
\min\limits_{[y^1;\ldots;y^K] \in Y_1\times \cdots\times Y_K}
\left\{\sum_{k=1}^K\left[\psi_k(y^k)+\Psi_k(y^k)\right]: \;\sum_{k=1}^K A_ky^k =b\right\}.
\end{equation}
Here for $1\leq k\leq K$ the domains $Y_k$ are closed and convex, $\psi_k(\cdot)$ are convex Lipschitz-continuous functions,
and $\Psi_k(\cdot)$ are {convex functions which are ``simple and fit'' $Y_k$.\footnote{The precise meaning of simplicity and fitting will be specified later. As of now, it suffices to give a couple of examples. When $\Psi_k$ is the $\ell_1$ norm, $Y_k$ can be the entire space, or the centered at the origin $\ell_p$-ball, $1\leq p\leq 2$; when $\Psi_k$ is the nuclear norm, $Y_k$ can be the entire space, or the centered at the origin Frobenius/nuclear norm ball.}}

\par
The problem of multi-term composite minimization (\ref{model:unconstrained:0}) has been considered (in a somewhat different setting) in \cite{Srebro12} for $K=2$. When $K=1$, problem (\ref{model:unconstrained:0}) becomes the usual composite minimization problem:
\begin{equation}\label{model:simplecomposition}
\min_{u\in U} \left\{\psi(u)+\Psi(u)\right\}
\end{equation}
which is well studied in the case where $\psi(\cdot)$ is a {\em smooth} convex function and $\Psi(\cdot)$ is a simple non-smooth function. For instance, it was shown that the composite versions of Fast Gradient Method originating in  Nesterov's seminal work \cite{NesCompMin} and further developed by many authors (see, e.g., \cite{Teboulle09a,R100,R101,Tseng08,Katya11} and references therein), as applied to (\ref{model:simplecomposition}), work {\sl as if} there were no nonsmooth term at all and exhibit the $O(1/t^2)$ convergence rate, which is the optimal rate attainable by first order algorithms of large-scale smooth convex optimization. Note that these algorithms cannot be directly applied to problems (\ref{model:unconstrained:0}) with $K>1$.

\par
The problem with semi-separable structures (\ref{model:constrained:0}) for $K=2$,  has also been extensively studied using the augmented Lagrangian approach (see, e.g., \cite{R102,Boyd10,Qin12,Goldfarb10,Goldfarb12,Goldfarb13,Monteiro13,Lan2014} and references therein). In particular, much work was carried out on the alternating directions method of multipliers (ADMM, see \cite{Boyd10} for an overview), which optimizes the augmented Lagrangian in an alternating fashion and exhibits an overall $O(1/t)$ convergence rate. Note that the available accuracy bounds for those algorithms involve optimal values of Lagrange multipliers of the equality constraints (cf. \cite{Lan2014}). Several variants of this method have been developed recently to adjust to the case for $K>2$ (see, e.g.\cite{Wo13}), however, most of these algorithms require to solve iteratively  subproblems of type (\ref{model:simplecomposition}) especially with the presence of non-smooth terms in the objective.

\subsection{{Our} contribution}
In this paper we do not assume smoothness of functions $\psi_k$, but, instead, we suppose that $\psi_k$ are {\em saddle point representable}:
\be
\psi_k(y^k)=\sup_{z^k\in Z_k}{[}\phi_k(y^k,z^k)-\oPsi_k(z^k){]},\;\;1\le k\le K,
\ee{model:saddlepoint1}
where $\phi_k(\cdot,\cdot)$ are smooth functions which are convex-concave (i.e., convex in their first and concave in the second argument), $Z_k$ are convex and compact, and $\oPsi_k(\cdot)$ are simple convex functions on $Z_k$. Let us consider, for instance, the multi-term composite minimization problem \rf{model:unconstrained:0}. Under \rf{model:saddlepoint1}, the primal problem \rf{model:unconstrained:0} allows for the saddle point reformulation:
\be
\min\limits_{y\in Y}\max\limits_{[z^1;\ldots;z^k]\in Z_1\times \cdots\times Z_K}\left\{
 \sum_{k=1}^K\left[\phi_k(A_ky+b_k,z^k)-\oPsi_k(z^k)+\Psi_k(A_ky+b_k)\right]\right\}
\ee{model:unconstrained}
Note that when there are no $\Psi_k,\oPsi_k$'s, problem \rf{model:unconstrained} becomes a {convex-concave} saddle point problem with smooth cost function, studied in \cite{FOM11}. In particular,
it was shown in \cite{FOM11} that Mirror Prox algorithm originating from \cite{MP}, when applied to the saddle point problem (\ref{model:unconstrained}), exhibits the  ``theoretically optimal'' convergence rate $O(1/t)$. Our goal in this paper is to develop novel $O(1/t)$-converging first order algorithms for problem \rf{model:unconstrained} (and also the related saddle point reformulation of the problem in \rf{model:constrained:0}), which appears to be the best rate known, under circumstances, from the literature (and established {there in essentially less general setting} than {the one} considered below).

\par
Our key observation is that composite problem {(\ref{model:unconstrained:0}), (\ref{model:saddlepoint1})} can be reformulated as a smooth linearly constrained saddle point problem by simply moving the nonsmooth terms into the problem domain. Namely, problem {(\ref{model:unconstrained:0}) , (\ref{model:saddlepoint1})} can be written as
\[
\min\limits_{{y\in Y,\; [y^k;\tau^k]\in Y_k^+\atop1\leq k\leq K}}\max\limits_{[z^k;\sigma^k]\in Z_k^+\atop1\leq k\leq K}\bigg\{\sum\limits_{k=1}^K\left[\phi_k(y^k,z^k)-\sigma^k+\tau^k\right]:
y^k=A_ky+b_k,\; k=1,...,K\bigg\}\\
\]
\[
Y_k^+=\left\{[y^k;\tau^k]: y^k\in Y_k, \tau^k\geq \Psi_k(y^k)\right\},Z_k^+=\left\{[z^k;\sigma^k]: z^k\in Z_k, \sigma^k\geq \oPsi_k(z^k)\right\},\;k=1,...,K.
\]
We can further approximate the resulting problem by penalizing
the equality constraints, thus passing to
\[
\min\limits_{{y\in Y,\; [y^k;\tau^k]\in Y_k^+\atop1\leq k\leq K}}\max\limits_{[z^k;\sigma^k]\in Z_k^+\atop1\leq k\leq K}\bigg\{\sum\limits_{k=1}^K\left[\phi_k(y^k,z^k)-\sigma^k+\tau^k +\rho_k\|y^k-A_ky-b_k\|_2\right]\bigg\}
\]
\begin{equation}\label{eqsaddle}
=\min\limits_{{y\in Y,\; [y^k;\tau^k]\in Y_k^+\atop1\leq k\leq K}}\max\limits_{w^k\in W_k,\;[z^k;\sigma^k]\in Z_k^+\atop1\leq k\leq K}\bigg\{\sum\limits_{k=1}^K\left[\phi_k(y^k,z^k)-\sigma^k+\tau^k +\rho_k\langle y^k-A_ky-b_k,w^k\rangle\right]\bigg\},
\end{equation}
where $\rho_k>0$ are penalty parameters and $W_k=\{w^k:\|w^k\|_2\leq 1\}, k=1,\ldots,K$.

We solve the convex-concave saddle point problem \rf{eqsaddle} with {\sl smooth} cost function by $O(1/t)$-converging Mirror Prox algorithm. It is worth to mention that if the functions $\phi_k$, $\Psi_k$ are Lipschitz continuous on the domains $A_kY+b_k$,  and $\rho_k$ are selected properly, the saddle point problem is exactly equivalent to the problem of interest.
\par
The monotone operator $F$ associated with the saddle point problem in (\ref{eqsaddle}) has a special structure: the variables can be split into two blocks $u$ (all $y$-, $z$- and $w$-variables) and $v$ (all $\tau$- and $\sigma$-variables) in such a way that the induced partition of $F$ is $F=[F_u(u);F_v]$ with the $u$-component $F_u$ depending solely on $u$ and constant $v$-component $F_v$. We demonstrate below that in this case the basic Mirror Prox algorithm admits a ``composite'' version which works essentially ``as if'' there were no $v$-component at all. This composite version of  Mirror Prox will be the working horse of all {subsequent developments}.

The main body of this paper is organized as follows. In section \ref{sect:preliminaries} we present required background on variational inequalities with monotone operators and convex-concave saddle points. In section \ref{sect:MirrorProx} we present and justify the composite Mirror Prox algorithm. In sections \ref{MultiTerm} and \ref{AltDir}, we apply our
approach to problems (\ref{model:unconstrained:0}), (\ref{model:saddlepoint1}) and (\ref{model:constrained:0}), (\ref{model:saddlepoint1}). In section \ref{sect:illustration}, we illustrate our approach (including numerical results) as applied to the motivating examples. {All proofs missing in the main body of the paper are relegated to the appendix.}


\section{Preliminaries: Variational Inequalities and Accuracy Certificates}\label{sect:preliminaries}
\paragraph{Execution protocols and accuracy certificates.}
Let $X$ be a nonempty closed convex set in a Euclidean space $E$ and $F(x):X\to E$ be a vector field.
\par
Suppose that we process $(X,F)$ by an algorithm which generates
a sequence of search points $x_t\in X$, $t=1,2,...$, and computes the vectors $F(x_t)$, so that after $t$ steps we have at our disposal {\sl $t$-step execution protocol}
$\cI_t=\{x_\tau,F(x_\tau)\}_{\tau=1}^t$. By definition, an {\sl accuracy certificate} for this protocol is simply a collection $\lambda^t=\{\lambda^t_\tau\}_{\tau=1}^t$ of nonnegative reals summing up to 1.
We associate with the protocol $\cI_t$ and accuracy certificate $\lambda^t$ two quantities as follows:\\
\begin{itemize}
\item {\sl Approximate solution} $x^t(\cI_t,\lambda^t):=\sum_{\tau=1}^t \lambda^t_\tau x_\tau$, which is a point of $X$;
\item {\sl Resolution $\Res(X'\big|\cI_t,\lambda^t)$ on a subset $X'\neq\emptyset$} of $X$ given by
\begin{equation}\label{resolution}
\Res(X'\big|\cI_t,\lambda^t) = \sup\limits_{x\in X'}\sum_{\tau=1}^t\lambda^t_\tau\langle F(x_\tau),x_\tau-x\rangle.
\end{equation}
\end{itemize}
The role of those notions in the optimization context is explained next\footnote{our exposition follows \cite{NOR}.}.
\paragraph{Variational inequalities.} Assume that $F$ is {\em monotone}, i.e.,
\be\langle F(x)-F(y),x-y\rangle\ge 0, \;\;\forall x,y\in X
\ee{eq:near-monotone}
and let our goal be to approximate a weak solution
 to the variational inequality (v.i.) $\hbox{vi}(X,F)$ associated with $(X,F)$; weak solution is defined as a point $x_*\in X$ such that
\begin{equation}\label{vin}
\langle F(y),y-x_*\rangle \geq0\,\,\forall y\in X.
\end{equation}
A natural (in)accuracy measure  of a candidate weak solution $x\in X$ to $\hbox{vi}(X,F)$ is the {\sl dual gap function}
\begin{equation}\label{dualgap}
\epsilonvi(x\big|X,F) = \sup_{y\in X} \langle F(y),x-y\rangle
\end{equation}
This inaccuracy is a  convex nonnegative function which vanishes exactly at the set of weak solutions to the \vixf.
\begin{proposition}\label{prop1} For every $t$, every execution protocol $\cI_t=\{x_\tau\in X,F(x_\tau)\}_{\tau=1}^t$ and every
accuracy certificate $\lambda^t$ one has $x^t:=x^t(\cI_t,\lambda^t)\in X$. Besides this, {assuming $F$ monotone,} for  every closed convex set $X'\subset X$ such that $x^t\in X'$ one has
\begin{equation}\label{epsilonvi}
\epsilonvi(x^t\big|X',F)\leq \Res(X'\big|\cI_t,\lambda^t).
\end{equation}
\end{proposition}
{\bf Proof.} Indeed, $x^t$ is a convex combination of the points $x_\tau\in X$ with coefficients $\lambda^t_\tau$, whence $x^t\in X$. With $X'$ as in the premise of Proposition, we have
$$
\forall y\in X': \langle F(y),x^t-y\rangle =\sum_{\tau=1}^t\lambda^t_\tau\langle F(y),x_\tau-y\rangle \leq \sum_{\tau=1}^t \lambda^t_\tau \langle F(x_\tau),x_\tau-y\rangle
\leq\Res(X'\big|\cI_t,\lambda^t),
$$
where the {first $\leq$ is due to} monotonicity of $F$. \qed
\paragraph{Convex-concave saddle point problems.} Now let $X=X_1\times X_2$, where $X_i$ is a closed convex subset in Euclidean space $E_i$, $i=1,2$, and $E=E_1\times E_2$, and let
$\Phi(x^1,x^2):X_1\times X_2\to\bR$
be a locally Lipschitz continuous function which is convex in $x^1\in X_1$ and concave in $x^2\in X_2$. $X_1,X_2,\Phi$ give rise to  the saddle point problem
\begin{equation}\label{saddlepoint}
\SadVal=\min_{x^1\in X_1}\max_{x^2\in X_2}\Phi(x^1,x^2),
\end{equation} two induced convex optimization problems
\begin{equation}\label{induced}
\begin{array}{rclr}
\Opt(P)&=&\min_{x^1\in X_1}\left[\overline{\Phi}(x^1)=\sup_{x^2\in X_2}\Phi(x^1,x^2)\right]&(P)\\
\Opt(D)&=&\max_{x^2\in X_2}\left[\underline{\Phi}(x^2)=\inf_{x^1\in X_1}\Phi(x^1,x^2)\right]&(D)\\
\end{array}
\end{equation}
and a vector field $F(x^1,x^2)=[F_1(x^1,x^2);F_2(x^1,x^2)]$ specified (in general, non-uniquely) by the relations
$$
\forall (x^1,x^2)\in X_1\times X_2: F_1(x^1,x^2)\in\partial_{x^1}\Phi(x^1,x^2),\,F_2(x^1,x^2)\in\partial_{x^2}[-\Phi(x^1,x^2)].
$$
It is well known that $F$ is monotone on $X$, and that weak solutions to the \vixf are exactly the saddle points of $\Phi$ on $X_1\times X_2$. These saddle points exist if and only if
$(P)$ and $(D)$ are solvable with equal optimal values, in which case the saddle points are exactly the pairs $(x^1_*,x^2_*)$ comprised by optimal solutions to $(P)$ and $(D)$.
In general, $\Opt(P)\geq\Opt(D)$, with equality definitely taking place when at least one of the sets $X_1,X_2$ is bounded; if both are bounded, saddle points do exist. {To avoid unnecessary complications, from now on, when speaking about a convex-concave saddle point problem, we assume that the problem is {\sl proper}, meaning that $\Opt(P)$ and $\Opt(D)$ are reals; this definitely is the case when $X$ is bounded.
\par
A natural (in)accuracy measure for a candidate $x=[x^1;x^2]\in X_1\times X_2$ to the role of a saddle point of $\Phi$ is the quantity
\begin{equation}\label{saddlepointinacc}
\begin{array}{rcl}
\epsilonsad(x\big|X_1,X_2,\Phi)&=&\overline{\Phi}(x^1) -\underline{\Phi}(x^2)\\
&=&[\overline{\Phi}(x^1)-\Opt(P)]+[\Opt(D)-\underline{\Phi}(x^2)] +\underbrace{[\Opt(P)-\Opt(D)]}_{\geq0}\\
\end{array}
\end{equation}
This inaccuracy is nonnegative and is the sum of the duality gap $\Opt(P)-\Opt(D)$ (always nonnegative and vanishing when one of the sets $X_1,X_2$ is bounded)  and the inaccuracies,
in terms of respective objectives, of $x^1$ as a candidate solution to $(P)$ and $x^2$ as a candidate solution to $(D)$.
\par
The role of accuracy certificates in convex-concave saddle point problems stems from the following observation:
\begin{proposition} \label{prop2} Let $X_1,X_2$ be nonempty closed convex sets, $\Phi:X:=X_1\times X_2\to\bR$ be a locally Lipschitz continuous convex-concave function, and
$F$ be the associated monotone vector field on $X$.
\par Let $\cI_t=\{x_\tau=[x^1_\tau;x^2_\tau]\in X,{F}(x_\tau)\}_{\tau=1}^t$ be a $t$-step execution protocol associated with $(X,{F})$ and
$\lambda^t=\{\lambda^t_\tau\}_{\tau=1}^t$ be an associated
accuracy certificate. Then $x^t:=x^t(\cI_t,\lambda^t)=[x^{1,t};x^{2,t}] \in X$.
\par
Assume, further, that $X^\prime_1\subset X_1$ and $X^\prime_2\subset X_2$ are closed convex sets such that
\begin{equation}\label{condition}
x^t\in X^\prime:=X^\prime_1\times X^\prime_2.
\end{equation}
Then
\begin{equation}\label{saddlepointres}
\epsilonsad(x^t\big|X_1^\prime,X_2^\prime,\Phi)=\sup_{x^2\in X_2'}\Phi(x^{1,t},x^2)-\inf_{x^1\in X_1'}\Phi(x^1,x^{2,t})\leq\Res(X^\prime\big|\cI_t,\lambda^t).
\end{equation}
In addition, setting $\widetilde{\Phi}(x^1)=\sup_{x^2\in X_2^\prime}\Phi(x^1,x^2)$, for every $\bar{x}^1\in X^\prime_1$ we have
\begin{equation}\label{forfeasible}
\widetilde{\Phi}(x^{1,t})-\widetilde{\Phi}(\bar{x}^1)\leq\widetilde{\Phi}(x^{1,t})-\Phi(\bar{x}^1,x^{2,t})\leq\Res(\{\bar{x}^1\}\times X^\prime_2\big|\cI_t,\lambda^t).
\end{equation}
In particular, when the problem $\Opt=\min_{x^1\in X^\prime_1}\widetilde{\Phi}(x^1)$ is solvable with an optimal solution $x^1_*$, we have
\begin{equation}\label{ifsolvable}
\widetilde{\Phi}(x^{1,t})-\Opt\leq\Res(\{x^1_*\}\times X^\prime_2\big|\cI_t,\lambda^t).
\end{equation}
\end{proposition}
{\bf Proof.} The inclusion $x^t\in X$ is evident. For every set $Y\subset X$ we have
$$
\begin{array}{l}
\forall [p;q]\in Y:\\
\Res(Y\big|\cI_t,\lambda^t)\geq \sum_{\tau=1}^t\lambda^t_\tau\left[\langle F_1(x^1_\tau),x^1_\tau-p\rangle + \langle F_2(x^2_\tau),x^2_\tau-q\rangle \right]\\
\geq \sum_{\tau=1}^t\lambda^t_\tau\left[[\Phi(x^1_\tau,x^2_\tau)-\Phi(p,x^2_\tau)]+[\Phi(x^1_\tau,q)-\Phi(x^1_\tau,x^2_\tau)]\right]\\
\hbox{[by the origin of $F$ and since $\Phi$ is convex-concave]}\\
=\sum_{\tau=1}^t\lambda^t_\tau\left[\Phi(x^1_\tau,q)-\Phi(p,x^2_\tau)\right]
\geq \Phi(x^{1,t},q)-\Phi(p,x^{2,t})\\
\;~\;~\;~\;~\;~\;~\;~\;~\;~\hbox{[by origin of $x^t$ and since $\Phi$ is convex-concave]}\\
\end{array}
$$
Thus, for every $Y\subset X$ we have
\begin{equation}\label{conclude}
\sup_{[p;q]\in Y} \left[\Phi(x^{1,t},q)-\Phi(p,x^{2,t})\right]\leq \Res(Y\big|\cI_t,\lambda^t).
\end{equation}
Now assume that (\ref{condition}) takes place. Setting $Y=X':=X^\prime_1\times X^\prime_2$ and recalling what $\epsilonsad$ is, (\ref{conclude}) yields (\ref{saddlepointres}). With $Y=\{\bar{x}^1\}\times X^\prime_2$ (\ref{conclude}) yields the second inequality in
(\ref{forfeasible}); the first inequality in (\ref{forfeasible}) is evident due to $x^{2,t}\in X^\prime_2$. \qed


\section{Composite Mirror Prox Algorithm}\label{sect:MirrorProx}
\subsection{The situation}\label{situation}
Let $U$ be a nonempty closed convex domain in a Euclidean space $E_u$, $E_v$ be a Euclidean space, and $X$ be a nonempty closed convex domain in $E=E_u\times E_v$. We denote vectors from $E$ by $x=[u;v]$ with blocks $u,v$ belonging to $E_u$ and $E_v$, respectively.\par
 We assume that
\begin{enumerate}
\item[{\bf A1}:] $E_u$ is equipped with a norm $\|\cdot\|$, the conjugate norm being $\|\cdot\|_*$, and $U$ is equipped with a {\sl distance-generating function} (d.g.f.) $\omega(\cdot)$ (that is, with a continuously differentiable convex function $\omega(\cdot):U\to\bR$)  which is {\sl compatible} with $\|\cdot\|$, meaning that $\omega$ is strongly convex, modulus 1, w.r.t. $\|\cdot\|$.
    \par
    Note that d.g.f. $\omega$ defines the {\sl Bregman distance}
    \be
    V_u(w):=\omega(w)-\omega(u)-\langle\omega'(u),w-u\rangle \geq {1\over 2}\|w-u\|^2,\,u,w\in U,
    \ee{eq:bregman}
    where the concluding inequality follows from strong convexity, modulus 1, of the d.g.f. w.r.t. $\|\cdot\|$.\par
    {In the sequel, we refer to the pair $\|\cdot\|,\,\omega(\cdot)$ as to {\sl proximal setup} for $U$.}
\item[{\bf A2}:] the image $PX$ of $X$ under the projection $x=[u;v]\mapsto Px:=u$ is contained in $U$.
\item[{\bf A3}:] we are given a vector field $F(u,v):X\to E$ on $X$ of the special structure as follows:
$$
F(u,v)=[F_u(u);F_v],
$$
with $F_u(u)\in E_u$ and $F_v\in E_v$. Note that $F$ is independent of $v$.
\par
We assume also that
\begin{equation}\label{ML}
\forall u,u'\in U: \|F_u(u)-F_u(u')\|_*\leq L\|u-u'\|+M
\end{equation}
with some $L<\infty$, $M<\infty$.
\item[{\bf A4}:] the linear form $\langle F_v,v\rangle$ of $[u;v]\in E$ is bounded from below on $X$ and is {coercive on $X$ w.r.t. $v$: whenever $[u_t;v_t]\in X$, $t=1,2,...$ is a sequence such that $\{u_t\}_{t=1}^\infty$ is bounded and $\|v_t\|_2\to\infty$ as $t\to\infty$, we have $\langle F_v,v_t\rangle \to\infty$, $t\to\infty$.}
\end{enumerate}

\paragraph{Our goal} in this section is to show that {\sl in the situation in question, proximal type processing $F$} (say, $F$ is monotone on $X$, and we want to solve the variational inequality given by $F$ and $X$) {\sl can be implemented ``as if'' there were no $v$-components in the domain and in $F$.}
\par
A generic application we are aiming at is as follows. We want to solve a ``composite'' saddle point problem
\begin{equation}\label{compositeSP}
\SadVal=\min_{u_1\in U_1}\max_{u_2\in U_2} \left[\phi(u_1,u_2)+\Psi_1(u_1)-\Psi_2(u_2)\right],
\end{equation}
where
\begin{itemize}
\item $U_1\subset E_1$ and $U_2\subset E_2$ are nonempty closed convex sets in Euclidean spaces $E_1,E_2$
\item $\phi$ is a smooth (with Lipschitz continuous gradient) convex-concave function on $U_1\times U_2$
\item $\Psi_1:U_1\to\bR$ and $\Psi_2:U_2\to\bR$ are convex functions, perhaps nonsmooth, but ``fitting'' the domains $U_1$, $U_2$ in the following
sense: for $i=1,2$, we can equip $E_i$ with a norm $\|\cdot\|_{(i)}$, and $U_i$ - with a compatible with this norm d.g.f.
$\omega_i(\cdot)$ in such a way that optimization problems of the form
\begin{equation}\label{auxil}
\min_{u_i\in U_i}\left[\alpha \omega_i(u_i)+\beta \Psi_i(u_i) +\langle \xi,u_i\rangle\right]\qquad{[\alpha>0,\beta>0]}
\end{equation}
are easy to solve.
\end{itemize}
Our ultimate goal is to solve (\ref{compositeSP}) ``as if'' there were no (perhaps) nonsmooth terms $\Psi_i$.
With our approach, we intend to ``get rid'' of the nonsmooth terms by ``moving'' them into the description of problem's domains. To this end,
we act as follows:
\begin{itemize}
 \item For $i=1,2$, we set $X_i=\{x_i=[u_i;v_i]\in E_i\times \bR: u_i\in U_i,v_i\geq\Psi_i(u_i)\}$ and set\\
 $$
 \begin{array}{c}
U:=U_1\times U_2\subset E_u:=E_1\times E_2,  E_v=\bR^2,\\
 X=\{x=[u=[u_1;u_2];v=[v_1;v_2]]:u_i\in U_i,v_i\geq \Psi_i(u_i),i=1,2\}\subset E_u\times E_v,
 \end{array}
$$
thus ensuring that $PX\subset U$, where $P[u;v]=u$;
\item We rewrite the problem of interest  equivalently as
\begin{equation}\label{equivalently}
\SadVal=\min_{x^1=[u_1;v_1]\in X_1}\max_{x^2=[u_2;v_2]\in X_2}\left[\Phi(u_1,v_1;u_2,v_2)=\phi(u_1,u_2)+v_1-v_2\right]
\end{equation}
Note that $\Phi$ is convex-concave and smooth. The associated monotone operator is
$$
F(u=[u_1;u_2],v=[v_1;v_2])=\left[F_u(u)=[\nabla_{u_1}\phi(u_1,u_2);-\nabla_{u_2}\phi(u_1,u_2)];F_v=[1;1]\right]
$$
and is of the structure required in {\bf A3}. Note that $F$ is Lipschitz continuous, so that (\ref{ML}) is satisfied with properly selected
$L$ and with $M=0$.
\end{itemize}
We intend to process the reformulated saddle point problem (\ref{equivalently}) with a properly modified state-of-the-art Mirror Prox (MP) saddle point
algorithm \cite{MP}. In its basic version and as applied to a variational inequality with Lipschitz continuous monotone operator (in particular, to a convex-concave
saddle point problem with smooth cost function), this algorithm exhibits $O(1/t)$ rate of convergence, which is the best rate achievable with First Order
saddle point algorithms as applied to large-scale saddle point problems (even those with bilinear cost function). The basic MP would require to
equip the domain $X=X_1\times X_2$ of (\ref{equivalently}) with a d.g.f. $\omega(x_1,x_2)$ resulting in an easy-to-solve auxiliary problems
of the form
\begin{equation}\label{auxil1}
\min_{x=[u_1;u_2;v_1;v_2]\in X}\left[ \omega(x) +\langle \xi,x\rangle\right],
\end{equation}
which would require to ``account nonlinearly'' for the $v$-variables (since $\omega$ should be a strongly convex  in both $u$-
and $v$-variables). While it is easy to construct $\omega$ from our postulated ``building blocks'' $\omega_1$, $\omega_2$ leading to easy-to-solve
problems (\ref{auxil}), this construction results in auxiliary problems (\ref{auxil1}) somehow more complicated than problems (\ref{auxil}).
To overcome this difficulty, below we develop a ``composite'' Mirror Prox algorithm taking advantage of the special structure of $F$, as expressed in {\bf A3},
and preserving the favorable efficiency estimates of the prototype. The modified MP operates with the auxiliary problems of the form
$$
\min_{x=[u_1;u_2;v_1;v_2]\in X_1\times X_2}\sum_{i=1}^2\left[\alpha_i \omega_i(u_i)+\beta_i v_i+\langle\xi_i,u_i\rangle\right],
\eqno{[\alpha_i>0,\beta_i>0]}
$$
that is, with pairs of uncoupled problems
$$
\min_{[u_i;v_i]\in X_i} \left[\alpha_i \omega_i(u_i)+\beta_i v_i+\langle\xi_i,u_i\rangle\right],\,i=1,2;
$$
recalling that $X_i=\{[u_i;v_i]:u_i\in U_i,v_i\geq\Psi_i(u_i)\}$, these problems are nothing but the easy-to-solve problems (\ref{auxil}).

\subsection{Composite Mirror-Prox algorithm}\label{MP}
Given the situation described in section \ref{situation}, we define the  associated {\em prox-mapping}: for $\xi=[\eta;\zeta]\in E$ and $x=[u;v]\in X$,
\be
P_x(\xi)\in \Argmin_{[s;w]\in X} \left\{\langle \eta-\omega'(u),s\rangle+\langle\zeta,w\rangle +\omega(s)\right\}\equiv\Argmin_{[s;w]\in X} \left\{\langle \eta,s\rangle+\langle\zeta,w\rangle +V_u(s)\right\}
\ee{eq:prox}
Observe that {$P_x([\eta;\gamma F_v])$ is well defined whenever $\gamma>0$} -- the {required} $\Argmin$ is nonempty due to the strong convexity of $\omega$ on $U$ and assumption {\bf A4} {(for verification, see item 0$^o$  in Appendix \ref{sect:prooftheMP}).}
Now consider the process as follows:
\begin{equation}\label{MirrorProx}
\begin{array}{rcl}
x_1:=[u_1;v_1]&\in& X;\\
y_{\tau}:=[u'_{\tau};v'_{\tau}]&=&P_{x_\tau}(\gamma_\tau F(x_\tau))=P_{x_\tau}(\gamma_\tau[ F_u(u_\tau);F_v])\\
x_{\tau+1}:=[u_{\tau+1};v_{\tau+1}]&=&P_{x_\tau}(\gamma_\tau F(y_\tau))=P_{x_\tau}(\gamma_\tau[ F_u(u'_\tau);F_v]),
\end{array}
\end{equation}
{where $\gamma_\tau>0$; the latter relation, due to the above, implies that the recurrence (\ref{MirrorProx}) is well defined.}

\begin{theorem}\label{theMP}
In the setting of section \ref{situation}, assuming that  {\bf A1}--{\bf A4} hold, consider the Composite Mirror Prox recurrence {\rm \rf{MirrorProx}} (CoMP) with stepsizes $\gamma_\tau>0$, $\tau=1,2,...$ satisfying the relation:
\be
\delta_\tau:=\gamma_\tau\langle F_u(u'_\tau)-F_u(u_\tau),u'_\tau-u_{\tau+1}\rangle-V_{u'_\tau}(u_{\tau+1})-V_{u_\tau}(u'_{\tau})\le \gamma_\tau^2M^2.
\ee{gammaupperbound}
Then the corresponding execution protocol $\cI_t=\{y_\tau,F(y_\tau)\}_{\tau=1}^t$
admits accuracy certificate $\lambda^t=\{\lambda^t_\tau=\gamma_\tau/\sum_{i=1}^t\gamma_i\}$ such that for every $X'\subset X$ it holds
\begin{equation}\label{itholdsXprime}
\Res(X'\big|\cI_t,\lambda^t)\leq {\Theta[X']+M^2\sum_{\tau=1}^t\gamma_\tau^2\over\sum_{\tau=1}^t\gamma_\tau},\;\;\;\Theta[X']=\sup_{[u;v]\in X'}V_{u_1}(u).
\end{equation}
Relation {\rm\rf{gammaupperbound}} is definitely satisfied when $0<\gamma_\tau\leq ({\sqrt{2}L})^{-1}$, or, in the case of $M=0$, when $\gamma_\tau\leq L^{-1}$.
\end{theorem}
Invoking Propositions \ref{prop1}, \ref{prop2}, we arrive at the following
\begin{corollary}\label{corMP}
Under the premise of Theorem \ref{theMP}, for every $t=1,2,...$, setting\\
$$
x^t=[u^t;v^t]={1\over\sum_{\tau=1}^t\gamma_\tau}\sum_{\tau=1}^t\gamma_\tau y_\tau.
$$
we ensure that  $x^t\in X$ and that\\
{\rm (i)} In the case  when $F$ is monotone on $X$, we have
\begin{equation}\label{eposilonvismall}
\epsilonvi(x^t\big|X,F)\leq \left[{\sum}_{\tau=1}^t\gamma_\tau\right]^{-1}\left[\Theta[X]+M^2{\sum}_{\tau=1}^t\gamma_\tau^2\right].
\end{equation}
{\rm (ii)} Let $X=X_1\times X_2$, and let $F$ be the monotone vector field associated with the saddle point problem {\rm (\ref{saddlepoint})}
with convex-concave locally Lipschitz continuous cost function $\Phi$. Then
\begin{equation}\label{eposilonsadsmall}
\epsilonsad(x^t\big|X_1,X_2,\Phi)\leq \left[{\sum}_{\tau=1}^t\gamma_\tau\right]^{-1}\left[\Theta[X]+M^2{\sum}_{\tau=1}^t\gamma_\tau^2\right].
\end{equation}
In addition, assuming  that problem $(P)$ in {\rm (\ref{induced})} is solvable with optimal solution $x^1_*$ and denoting by $x^{1,t}$ the projection of $x^t\in X=X_1\times X_2$ onto $X_1$, we have
\be
\overline{\Phi}(x^{1,t})-\Opt(P)\leq \left[{\sum}_{\tau=1}^t\gamma_\tau\right]^{-1}\left[\Theta[\{x^1_*\}\times X_2]+M^2{\sum}_{\tau=1}^t\gamma_\tau^2\right].
\ee{optimalitygap}
\end{corollary}
\begin{remark}\label{rem1} When $F$ is Lipschitz continuous (that is, {\rm (\ref{ML})} holds true with some $L>0$ and $M=0$), the requirements on the
stepsizes imposed in the premise of Theorem \ref{theMP}
reduce to $\delta_\tau\leq0$ for all $\tau$ and are definitely satisfied with the constant stepsizes $\gamma_\tau=1/L$.
Thus, in the case under consideration we can assume w.l.o.g. that $\gamma_\tau\geq 1/L$, thus ensuring that the upper bound
on $\Res(X'\big|\cI_t,\lambda^t)$ in {\rm (\ref{itholdsXprime})}
is $\leq \Theta[X']Lt^{-1}$. As a result, {\rm (\ref{optimalitygap})} becomes
\begin{equation}\label{optimalitygap1}
\overline{\Phi}(x^{1,t})-\Opt(P)\leq \Theta[\{x^1_*\}\times X_2]Lt^{-1}.
\end{equation}
\end{remark}

\subsection{Modifications}
In this section, we demonstrate that in fact our algorithm admits some freedom in building approximate solutions, freedom which can be used to improve to some extent solutions' quality. Modifications to be presented originate from \cite{NR02}.  We assume that we are in the situation described in section \ref{situation}, and assumptions {\bf A1} -- {\bf A4} are in force. In addition, we assume that
\begin{enumerate}
\item[{\bf A5}:] The vector field $F$ described in {\bf A3} is monotone, and the variational inequality given by $(X,F)$ has a weak solution:
\begin{equation}\label{weaksol}
\exists x_*=[u_*;v_*]\in X: \langle F(y),y-x_*\rangle \geq 0\,\,\forall y\in X
\end{equation}
\end{enumerate}

\begin{lemma}\label{lemnew1}
In the situation from section \ref{situation} and under assumptions {\bf A1} -- {\bf A5}, for $R\geq0$ let us set
\begin{equation}\label{widehatTheta}
\widehat{\Theta}(R)=\max_{u,u'\in U}\left\{V_u(u'):\|u-u_1\|\leq R, \|u'-u_1\|\leq R\right\}
\end{equation}
(this quantity is finite since $\omega$ is continuously differentiable on $U$), and let \[
\{x_\tau=[u_\tau;v_\tau]:\tau\leq N+1,y_\tau:\tau\leq N\}\]
 be the trajectory of the $N$-step Mirror Prox algorithm {\rm (\ref{MirrorProx})} with stepsizes $\gamma_\tau>0$ which ensure {\rm (\ref{gammaupperbound})} for $\tau\leq N$. Then for all
 $u\in U$ and $t\leq N+1$,
\begin{equation}\label{thenwehave}
 0\leq V_{u_t}(u)\leq \widehat{\Theta}(\max[R_N,\|u-u_1\|]),\;\;\;R_N:=2\left(2V_{u_1}(u_*)+M^2\sum_{\tau=1}^{N-1}
 \gamma_\tau^2\right)^{1/2},
\end{equation}
with $u_*$ defined in {\rm (\ref{weaksol})}.
\end{lemma}

\begin{proposition}\label{propDifferentWeights}
In the situation of section \ref{situation} and under assumptions {\bf A1} -- {\bf A5}, let $N$ be a positive integer, and let  $\cI_N=\{y_\tau,F(y_\tau)\}_{\tau=1}^N$ be the execution protocol generated by $N$-step Composite Mirror Prox recurrence {\rm (\ref{MirrorProx})} with stepsizes
$\gamma_\tau$ ensuring  {\rm (\ref{gammaupperbound})}. Let also $\lambda^N=\{\lambda_1,...,\lambda_N\}$ be a collection of positive reals summing up to 1 and such that
\begin{equation}\label{monotone}
\lambda_1/\gamma_1\leq\lambda_2/\gamma_2\leq...\leq\lambda_N/\gamma_N.
\end{equation}
Then for every $R\geq0$, with $X_R=\{x=[u;v]\in X: \|u-u_1\|\leq R\}$ one has
\begin{equation}\label{newbound}
\Res(X_R|\cI_N,\lambda^N)\leq {\lambda_N\over\gamma_N}\widehat{\Theta}(\max[R_N,R])+M^2\sum_{\tau=1}^N\lambda_\tau\gamma_\tau,
\end{equation}
with $\widehat{\Theta}(\cdot)$ and $R_N$ defined by {\rm (\ref{widehatTheta}) and (\ref{thenwehave})}.
\end{proposition}

Invoking Propositions \ref{prop1}, \ref{prop2}, we arrive at the following modification of Corollary \ref{corMP}.
\begin{corollary}\label{corMPNew} Under the premise and in the notation of Proposition \ref{propDifferentWeights}, setting\\
$$
x^N=[u^N;v^N]=\sum_{\tau=1}^N\lambda_\tau y_\tau.
$$
we ensure that  $x^N\in X$. Besides this,  \\
{\rm (i)} Let $X'$ be a closed convex subset of $X$ such that $x^N\in X'$ and the projection of $X'$ on the $u$-space is contained in $\|\cdot\|$-ball of radius $R$ centered at $u_1$. Then
\begin{equation}\label{eposilonvismallNew}
\epsilonvi(x^N\big|X',F)\leq {\lambda_N\over\gamma_N}\widehat{\Theta}(\max[R_N,R])+M^2\sum_{\tau=1}^N\lambda_\tau\gamma_\tau.
\end{equation}
{\rm (ii)} Let $X=X_1\times X_2$ and $F$ be the monotone vector field associated with saddle point problem {\rm (\ref{saddlepoint})}
with convex-concave locally Lipschitz continuous cost function $\Phi$. Let, further, $X^\prime_i$ be closed convex subsets of $X_i$, $i=1,2$, such that $x^N\in X^\prime_1\times X^\prime_2$ and the projection of $X^\prime_1\times X^\prime_2$ onto the $u$-space is contained in $\|\cdot\|$-ball of radius $R$ centered at $u_1$. Then
\begin{equation}\label{eposilonsadsmallNew}
\epsilonsad(x^N\big|X^\prime_1,X^\prime_2,\Phi)\leq {\lambda_N\over\gamma_N}\widehat{\Theta}(\max[R_N,R])+M^2{\sum}_{\tau=1}^N\lambda_\tau\gamma_\tau.
\end{equation}
\end{corollary}

\section{Multi-Term Composite Minimization}\label{MultiTerm}
{In this section, we focus on the problem  (\ref{model:unconstrained:0}), (\ref{model:saddlepoint1}) of multi-term composite minimization.}
\subsection{Problem setting}\label{sec:41}
{We intend to consider problem (\ref{model:unconstrained:0}), (\ref{model:saddlepoint1}) in the situation as follows. For a nonnegative integer $K$ and $0\le k\le K$ we are given
\begin{enumerate}
\item 
Euclidean spaces $E_k$ and $\oE_k$ along with their nonempty closed convex subsets $Y_k$ and $Z_k$, respectively;
\item 
Proximal setups for $(E_k,Y_k)$ and $(\oE_k, Z_k)$, that is, norms $p_k(\cdot)$ on $E_k$, norms $q_k(\cdot)$ on $\oE_k$, and {d.g.f.'s} $\omega_k(\cdot):Y_k\to \bR$, $\oomega_k(\cdot):Z_k\to \bR$ {compatible with $p_k(\cdot)$ and $q_k(\cdot)$, respectively;}
\item 
Affine mappings $y^0\mapsto A_k y^0+b_k:E_0\to E_k$, where $y^0\mapsto A_0y^0+b_0$ is the identity mapping on $E_0$;
\item 
Lipschitz continuous convex functions $\psi_k(y^k):Y_k\to\bR$ along with their \emph{saddle point representations}
\be
\psi_k(y^k)=\sup_{z^k\in Z_k}{[}\phi_k(y^k,z^k)-\oPsi_k(z^k){]},\;\;0\le k\le K,
\ee{model:saddlepoint2}
where $\phi_k(y^k,z^k):Y_k\times Z_k\to\bR$ are smooth (with Lipschitz continuous gradients)  functions convex in $y^k\in Y_k$ and concave in $z^k\in Z_k$, and  $\oPsi_k(z^k):Z_k\to\bR$ are Lipschitz continuous convex functions such that the problems of the form
\begin{equation}
\begin{array}{rl}\min\limits_{z^k\in Z_k}&\left[\oomega_k(z^k)+\langle \xi^k,z^k\rangle +\alpha\oPsi_k(z^k)\right]\quad[\alpha>0]
\end{array}
\end{equation}
are easy to solve;
\item 
Lipschitz continuous convex functions $\Psi_k(y^k):Y_k\to\bR$ such that the problems of the form
\begin{equation}\label{octeq20}
\begin{array}{rl}\min\limits_{y^k\in Y_k}&\left[\omega_k(y^k)+\langle \xi^k,y^k\rangle +\alpha\Psi_k(y^k)\right]\quad[\alpha>0]
\end{array}
\end{equation}
are easy to solve;
\item For $1\leq k\leq K$, the norms $\pi_k^*(\cdot)$ on $E_k$ are given, with conjugate norms $\pi_k(\cdot)$, along with {d.g.f.'s} $\widehat{\omega}_k(\cdot):\;W_k:=\{w^k\in E_k:\pi_k(w^k)\leq1\}\to\bR$ which are strongly convex, modulus 1, w.r.t. $\pi_k(\cdot)$ such that the problems
\begin{equation}\label{octeq21}
\min_{w^k\in W_k}\left[\widehat{\omega}_k(w^k)+\langle \xi^k,w^k\rangle\right]
\end{equation}
are easy to solve.
\end{enumerate}
The outlined data define the sets
\[
\begin{array}{rcl}
Y^+_k&=&\{[y^k;\tau^k]: \;y^k\in Y_k,\tau^k\geq\Psi_k(y^k)\}\subset E_k^+:=E_k\times \bR,\,\,0\leq k\leq K,\\
Z^+_k&=&\{[z^k;\sigma^k]: \;z^k\in Z_k,\sigma^k\geq\oPsi_k(z^k)\}\subset \oE_k^+:=\oE_k\times \bR,\,\,0\leq k\leq K.
\end{array}
\]

\paragraph{The problem of interest} {(\ref{model:unconstrained:0}), (\ref{model:saddlepoint1}) along with its saddle point reformulation in the just defined situation read}
\begin{subequations}
\begin{align}
\Opt&=\min\limits_{y^0\in Y_0}\left\{ f(y^0) : =
\sum\limits_{k=0}^K\left[\psi_k(A_ky^0+b_k)+\Psi_k(A_ky^0+b_k)\right]\right\} \label{octeq23a}\\
&=\min\limits_{y^0\in Y_0}\left\{ f(y^0) =\max\limits_{\{z^k\in Z_k\}_{k=0}^K}\sum\limits_{k=0}^K\left[\phi_k(A_ky^0+b_k,z^k)+\Psi_k(A_ky^0+b_k)-\oPsi_k(z^k)\right]\right\}\label{octeq23b}
\end{align}
\text{which we rewrite equivalently as}
\begin{align}\Opt
=\min\limits_{\{[y^k;\tau^k]\}_{k=0}^K\atop \in Y_0^+\times\cdots\times Y_K^+}\max\limits_{\{[z^k;\sigma^k]\}_{k=0}^K \atop \in Z_0^+\times\cdots\times Z_K^+}
\left\{\sum\limits_{k=0}^K\left[\phi_k(y^k,z^k)+\tau^k-\sigma^k\right]: y^k=A_ky^0+b_k,\;1\leq k\leq K\right\}.\label{octeq23c}
\end{align}
\end{subequations}

From now on we make the following assumptions
\begin{quote}
 {\bf B1}: {\sl We have $A_kY_0+b_k\subset Y_k$, $1\leq k\leq K$;}\\
{\bf B2}: {\sl For $0\leq k\leq K$, the sets $Z_k$ are bounded. Further, {the functions $\Psi_k$ are below bounded on $Y_k$, and} the functions {$f_k=\psi_k+\Psi_k$ }  { are coercive on} $Y_k$: whenever $y^k_t\in Y_k$, $t=1,2,...$, are such that $p_k(y^k_t)\to\infty$ as $t\to\infty$, we have
${f_k}(y^k_t)\to\infty$}.
\end{quote}
Note that  {\bf B1} and {\bf B2} imply that the saddle point problem \rf{octeq23c} is solvable; let  $\{[y^k_*;\tau^k_*]\}_{0\leq k\leq K};\;\{[z^k_*;\sigma^k_*]\}_{0\leq k\leq K} $ be the corresponding saddle point.

\subsection{Course of actions} Given $\rho_k>0$, $1\leq k\leq K$, we approximate (\ref{octeq23c}) by the problem

\begin{subequations}
\begin{align}
\widehat{\Opt}&=\min\limits_{\{[y^k;\tau^k]\}_{k=0}^K\atop \in Y_0^+\times\cdots\times Y_K^+}\max\limits_{\{[z^k;\sigma^k]\}_{k=0}^K \atop \in Z_0^+\times\cdots\times Z_K^+}
\left\{\sum\limits_{k=0}^K\left[\phi_k(y^k,z^k)+\tau^k-\sigma^k\right]+\sum\limits_{k=1}^K
\rho_k\pi_k^*(y^k-A_ky^0)\right\}. \label{octeq24a}\\
&=\min\limits_{x^1\in X_1 \atop:=Y_0^+\times\cdots\times Y_K^+}\max\limits_{x^2\in  X_2 \atop:= Z_0^+\times\cdots\times Z_K^+\times W_1\times\cdots W_K}\Phi\bigg(\underbrace{\{[y^k;\tau^k]\}_{k=0}^K}_{x^1},\;\underbrace{\left[\{[z^k;\sigma^k]\}_{k=0}^K; \{w^k\}_{k=1}^K\right]}_{x^2}\bigg)\label{octeq24b}
\end{align}
\end{subequations}
where
\[
\Phi(x^1,x^2)=\sum\limits_{k=0}^K\left[\phi_k(y^k,z^k)+\tau^k-\sigma^k\right]
+\sum\limits_{k=1}^K
\rho_k\langle w^k,y^k-A_ky^0-b_k\rangle.
\]
Observe that the monotone operator $F(x^1,x^2)=[F_1(x^1,x^2);F_2(x^1,x^2)]$ associated with the saddle point problem in (\ref{octeq24b}) is given by
\begin{equation}\label{octeq25}
\begin{array}{rcl}
F_1(x^1,x^2)&=&\bigg[\nabla_{y^0}\phi_0(y^0,z^0)-\sum\limits_{k=1}^K\rho_kA_k^Tw^k;1;\;
\left\{\nabla_{y^k}\phi_k(y^k,z^k)+\rho_kw^k;1\right\}_{k=1}^K\bigg],\\
F_2(x^1,x^2)&=&\bigg[\left\{-\nabla_{z^k}\phi_k(y^k,z^k);1\right\}_{k=0}^K;\;
\left\{-\rho_k[y^k-A_ky^0-b_k]\right\}_{k=1}^K\bigg].
\end{array}
\end{equation}
Now let us set
\begin{itemize}
\item $U=\left\{\begin{array}{rl}u=[y^0;...;y^K;z^0;...;z^K;w^1;...;w^K]:&y^k\in Y_k, \,z^k\in Z_k,\,0\leq k\leq K,\\
    &\pi_k(w^k)\leq1,\,1\leq k\leq K\end{array}\right\}$
\item $X=\left\{\begin{array}{rl} x=\left[u=[y^0;...;y^K;z^1;...;z^K;w^1;...;w^K];\;v=[\tau^0;...;\tau^K;\sigma^0;...;\sigma^K]\right]:\\ u\in U, \,\tau^k\geq \Psi_k(y^k),\,\sigma^k\geq \oPsi_k(z^k),\,0\leq k\leq K\end{array}\right\}$,\\
so that $PX\subset U$, cf. assumption {\bf A2} in section \ref{situation}.
\end{itemize}
The variational inequality associated with the saddle point problem in (\ref{octeq24b}) can be treated as the variational inequality on the domain $X$ with the monotone operator
\[F(x=[u;v])=[F_u(u);F_v],
\] where
\begin{equation}\label{octeq26}
\begin{array}{rl}
F_u(\underbrace{[y^0;...;y^K;\,z^0;...;z^K;\,w^1;...;w^{K}]}_{u})
=&\left[
\begin{array}{l}\nabla_y\phi_0(y^0,z^0)-\sum\limits_{k=1}^{K}\rho_kA_k^Tw^k\\
\left\{\nabla_y\phi_k(y^k,z^k)+\rho_kw^k\right\}_{k=1}^K\\
\left\{-\nabla_z\phi_k(y^k,z^k\right\}_{k=0}^K\\
\left\{-\rho_k[y^k-A_ky^0-b_k]\right\}_{k=1}^K\end{array}\right]\\
F_v(\underbrace{[\tau^0;...;\tau^K;\sigma^0;...;\sigma^K]}_{v}) =&[1;...;1].\\
\end{array}
\end{equation}
This operator meets the structural assumptions {\bf A3} and {\bf A4} from section \ref{situation} ({\bf A4} is guaranteed by {\bf B2}).
We can equip $U$ and its embedding space $E_u$ with the proximal setup $\|\cdot\|,\;\omega(\cdot)$ given by
\begin{equation}\label{aggregation}
\begin{array}{rcl}
\|u
\|&=&\sqrt{\sum_{k=0}^K[\alpha_kp_k^2(y^k)
+\beta_kq_k^2(z^k)]+\sum_{k=1}^{K}\gamma_k\pi_k^2(w^k)},\\
\omega(u
)&=&\sum_{k=0}^K[\alpha_k\omega_k(y^k)+\beta_k\oomega_k(z^k)]+\sum_{k=1}^{K}\gamma_k\widehat{\omega}_k(w^k),\\
\end{array}
\end{equation}
where $\alpha_k,\beta_k$, $0\leq k\leq K$, and  $\gamma_k$, $1\leq k\leq K$, are positive aggregation parameters.  Observe that carrying out a step of the CoMP
algorithm presented in  section \ref{MP} requires computing $F$ at $O(1)$ points of $X$  and solving $O(1)$ auxiliary problems of the form
\bse
\lefteqn{\min\limits_{{[y^0;...;y^K;z^0;...;z^K],\atop[;w^1;...;w^K;\tau^0;...;\tau^K;\sigma^0;...\sigma^K]}}
\left\{
\sum\limits_{k=0}^K \left[a_k\omega_k(y^k) +\langle\xi_k,y^k\rangle +b_k\tau^k\right]\right.}\\
&\left.+\sum\limits_{k=0}^K \left[c_k\oomega_k(z^k) +\langle\eta_k,z^k\rangle +d_k\sigma^k\right]+\sum\limits_{k=1}^{K}\left[
e_k\widehat{\omega}_k(w^k)+\langle \zeta_k,w^k\rangle\right]\right\}:\\
&y^k\in Y_k,\tau^k\geq \Psi_k(y^k),\;z^k\in Z_k,\sigma^k\geq \oPsi_k(y^k),\;0\leq k\leq K,\; \pi_k(w^k)\leq1,\,1\leq k\leq K,
\ese
with positive $a_k,...,e_k$, and we have assumed that these problems are easy to solve.

\subsection{``Exact penalty''}\label{ExactPenalty} Let us make one more assumption:
\begin{quote}
{\bf C}: {\sl For $1\leq k\leq K$, \begin{itemize}
\item $\psi_{k}$ are Lipschitz continuous on $Y_k$  with constants $G_k$ w.r.t. $\pi_k^*(\cdot)$,
\item $\Psi_{k}$ are Lipschitz continuous on $Y_k$ with constants $H_k$ w.r.t. $\pi_k^*(\cdot)$.
 \end{itemize}}
 \end{quote}
 Given a feasible solution $\overline{x}=[\overline{x}^1; \overline{x}^2]$,
$\overline{x}^1:=\{[\overline{y}^k;\overline{\tau}^k]\in Y^+_k\}_{k=0}^K$ to the saddle point problem \rf{octeq24b},
 let us set
$$
\widehat{y}^0=\overline{y}^0;\;\widehat{y}^k=A_k\overline{y}^0+b_k,\;1\leq k\leq K;\;\widehat{\tau}^k=\Psi_k(\widehat{y}^k),\;0\leq k\leq K,
 $$
thus getting another feasible (by assumption  {\bf B1}) solution $\widehat{x}=\big[\widehat{x}^1=\{[\widehat{y}^k;\widehat{\tau}^k]\}_{k=0}^K;\,\overline{x}^2\big]$  to (\ref{octeq24b}).  We call $\widehat{x}^1$ {\sl correction} of $\overline{x}^1$. For $1\leq k\leq K$ we  clearly have
\bse
\psi_k(\widehat{y}^k)&\leq& \psi_k(\overline{y}^k)+G_k\pi_k^*(\widehat{y}^k-\overline{y}^k)=\psi_k(\overline{y}^k)+G_k\pi_k^*(\overline{y}^k-A_k\overline{y}^0-b_k),\\
\widehat{\tau}^k&=&\Psi_k(\widehat{y}^k)\leq \Psi_k(\overline{y}^k)+H_k\pi_k^*(\widehat{y}^k-\overline{y}^k)\leq \overline{\tau}^k+H_k\pi_k^*(\overline{y}^k-A_k\overline{y}^0-b_k),
\ese
and $\widehat{\tau}^0=\Psi_0(\overline{y}^0)\leq\overline{\tau}^0$. Hence for
$
\overline{\Phi}(x^1)=\max\limits_{x^2\in X_2} \Phi(x^1,x^2)$ we have
$$
\overline{\Phi}(\widehat{x}^1)\leq\overline{\Phi}(\overline{x}^1)+\sum_{k=1}^{K}[H_k+G_k]\pi_k^*(\overline{y}^{k}-A_k\overline{y}^0-b_k)
-\sum_{k=1}^{K}\rho_k\pi_k^*(\overline{y}^{k}-A_k\overline{y}^0-b_k).
$$
We see that under the condition
\begin{equation}\label{octeq28}
\rho_k\geq G_k+H_k,\,\, 1\leq k\leq K,
\end{equation}
correction does not increase the value of the primal objective of (\ref{octeq24b}), whence the saddle point value $\widehat{\Opt}$ of (\ref{octeq24b}) is $\geq$ the optimal value $\Opt$ in the problem of interest
(\ref{octeq23a}). Since the opposite inequality is evident, we arrive at the following
\begin{proposition}\label{pro:correction}{In the situation of section \ref{sec:41} let assumptions {\bf B1}, {\bf B2}, {\bf C} and {\rm (\ref{octeq23a})} hold true.} Then

\begin{enumerate}[(i)]
\item  the optimal value $\widehat{\Opt}$ in {\rm (\ref{octeq24a})} coincides with the optimal value $\Opt$ in the problem of interest {\rm (\ref{octeq23a})};

\item  consequently, if $\overline{x}=[\overline{x}^1;\overline{x}^2]$ is a feasible solution of the saddle point problem in {\rm (\ref{octeq24b})}, then the correction
$\widehat{x}^1=\{[\widehat{y}^k;\widehat{\tau}^k]\}_{k=0}^K$ of  $\overline{x}^1$ is a feasible solution to the problem of interest {\rm (\ref{octeq23c})}, and
\be
f(\widehat{y}^0)-\Opt\leq \epsilonsad(\overline{x}\big|X_1,X_2,\Phi),
\ee{thebetter}
where $\widehat{y}^0(=y^0(\widehat{x}^1))$ is the ``$y^0$-component'' of $\widehat{x}^1$;
\end{enumerate}
\end{proposition}

As a corollary, {under the premise of Proposition \ref{pro:correction}}, when applying to the saddle point problem (\ref{octeq24b}) the CoMP algorithm induced by the above setup and passing ``at no cost'' from the approximate solutions $x^t=[x^{1,t};x^{2,t}]$ generated by CoMP to the corrections $\widehat{x}^{1,t}$ of $x^{1,t}$'s, we get feasible solutions to the problem of interest
(\ref{octeq23a}) satisfying the error bound
\begin{equation}\label{errorbound}
f({y}^0(\widehat{x}^{1,t}))-\Opt\leq {\Theta[x_*^1\times X_2]L\over t},\,t=1,2,...
\end{equation}
where $L$ is the Lipschitz constant of $F_u(\cdot)$ induced by the norm $\|\cdot\|$ given by (\ref{aggregation}), and $\Theta[\cdot]$ is induced by the {d.g.f.} given by the same (\ref{aggregation}) and the $u=[y^0;...;y^K;z^0;...z^K;w^1;...;w^K]$ -component of the starting point. Note that $W_k$ and $Z_k$ are compact, whence  $\Theta[x_*^1\times X_2]$ is finite.
\paragraph{Remark.} Note that the value of the penalty in \rf{octeq28} which guarantees the validity of correction (the bound \rf{thebetter} of Proposition \ref{pro:correction}) may be very conservative. When implementing the algorithm the coefficients $\rho_k$ of penalization can be adjusted on-line. Indeed, let $\overline{\Phi}(\overline{x}^1)=\sup_{x^2\in X_2}\Phi(\overline{x}^1,x^2)$ (cf \rf{induced}). We always have $\widehat{\Opt}\leq\Opt$.  It follows that independently of how $\rho_k$ are selected, we have
    \begin{equation}\label{epsilonone}
    f(\widehat{y}^0)-\Opt\leq \underbrace{[f(\widehat{y}^0)-\overline{\Phi}(\overline{x}^1)]}_{\epsilon_1}+
    \underbrace{\left[\overline{\Phi}(\overline{x}^1)-\widehat{\Opt}\right]}_{\epsilon_2}
    \end{equation}
    for every feasible solution $\overline{x}^1=\{[\overline{y}^k;\overline{\tau}^k]\}_{k=0}^K$ to  (\ref{octeq24b}) and the same inequality holds for its correction $\widehat{x}^1=\{[\widehat{y}^k;\widehat{\tau}^k]\}_{k=0}^K$.
    When $\overline{x}^1$ is a component of a good (with small $\epsilonsad$) approximate solution to the saddle point problem (\ref{octeq24b}), $\epsilon_2$ is small. If $\epsilon_1$ also is small, we are done; otherwise we can either increase in a fixed ratio the current values of all $\rho_k$, or only of those $\rho_k$ for which passing from $[\overline{y}^k;\overline{\tau}^k]$ to $[\widehat{y}^k;\widehat{\tau}^k]$ results in ``significant'' quantities
    $$
[\psi_k(\widehat{y}^k)+\widehat{\tau}^k]-[\psi_k(\overline{y}^k)+\overline{\tau}^k+\rho_k\pi_k^*(\overline{y}^k-A_k\overline{y}^0-b_k)]
$$
and solve the updated saddle point problem  (\ref{octeq24b}).

\subsection{Numerical illustrations}\label{sect:illustration}

\subsubsection{Matrix completion}\label{sec:LS}
\label{sec:mcomplet}
\paragraph{Problem of interest.} In the experiments to be reported, we applied the just outlined approach
to {Example 1, that is, to the problem}
\begin{equation}\label{twopenalties}
\begin{array}{c}
\Opt=\min\limits_{y^0\in\bR^{n\times n}} \big[\upsilon(y^0)=\underbrace{{1\over 2}\|P_\Omega y^0 -b\|_2^2}_{\psi_0(y^0)} +\underbrace{\lambda \|y^0\|_1}_{\Psi_0(y^0)}+\underbrace{\mu\|y^0\|_\nuc}_{\Psi_1(y^0)}\big].\\
\end{array}
\end{equation}
where $\Omega$ is a given set of cells in an $n\times n$ matrix, and $P_\Omega y$ is the restriction
of $y\in\bR^{n\times n}$ onto $\Omega$; this restriction is treated as a vector from $\bR^M$, $M=\Card(\Omega)$. Thus, (\ref{twopenalties}) is a kind of matrix completion problem where we want to recover a sparse and low rank $n\times n$ matrix given noisy observations $b$ of its entries in cells from $\Omega$.  Note that  (\ref{twopenalties}) is a special case of (\ref{octeq23b}) with $K=1$, $Y_0=Y_1=E_0=E_1=\bR^{n\times n}$, the identity mapping $y^0\mapsto A_1y^0$, and $\phi_0(y^0,z^0)\equiv \psi_0(y^0)$, $\phi_1\equiv 0$ (so that $Z_k$ can be defined as singletons, and $\oPsi_k(\cdot)$ set to 0, $k=0,1$).
\paragraph{Implementing the CoMP algorithm.}
When implementing the CoMP algorithm, we used the Frobenius norm $\|\cdot\|_F$ on $\bR^{n\times n}$ in the role of $p_0(\cdot)$, $p_1(\cdot)$ and $\pi_1(\cdot)$, and the function
${1\over 2}\|\cdot\|_F^2$ in the role of d.g.f.'s $\omega_0(\cdot)$, $\omega_1(\cdot)$, $\widehat{\omega}_1(\cdot)$.
\par
The aggregation weights in (\ref{aggregation}) were chosen as $\alpha_0=\alpha_1=1/D$ and $\gamma_1=1$, where $D$ is a {\sl guess} of the quantity $D_*:=\|y^0_*\|_F$, where $y^0_*$ is the optimal solution (\ref{twopenalties}). With $D=D_*$, our aggregation would roughly optimize the right hand side in (\ref{errorbound}), provided the starting point is the origin.
\par
The coefficient $\rho_1$ in (\ref{octeq24b}) was adjusted dynamically as explained at the end of section \ref{ExactPenalty}. Specifically, we start with a small (0.001) value of $\rho_1$ and restart the solution process, increasing by factor 3 the previous value of $\rho_1$, each time when  the $x^1$-component $\overline{x}$ of current approximate solution and its correction $\widehat{x}$ violate the inequality
$\upsilon(y^0(\widehat{x}))\leq (1+\kappa)\overline{\Phi}(\overline{x})$ for some small tolerance $\kappa$ (we used $\kappa=1.$e-4), cf. (\ref{epsilonone}).
\par
The stepsizes $\gamma_t$ in the CoMP algorithm were adjusted dynamically, specifically, as follows. At a step $\tau$, given a current guess $\gamma$ for the stepsize,
 we set $\gamma_\tau=\gamma$, perform the step and check whether $\delta_\tau\leq0$. If this is the case, we pass to step $\tau+1$,  the new guess for the stepsize being $1.2$ times the old one. If $\delta_\tau$ is positive, we decrease $\gamma_\tau$ in a fixed proportion (in our implementation -- by factor 0.8), repeat the step, and proceed in this fashion until the resulting value of $\delta_\tau$ becomes nonpositive. When it happens, we pass to step $\tau+1$, and use the value of $\gamma_\tau$ we have ended up with as our new guess for the stepsize.
 \par
 In all our experiments, the starting point was given by the matrix $\bar{y}:=P_\Omega^*b$ (``observations of entries in cells from $\Omega$ and zeros in all other cells'') according to $y^0=y^1=\bar{y}$, $\tau^0=\lambda\|\bar{y}\|_1$, $\tau^1=\mu\|\bar{y}\|_\nuc$, $w^1=0$.

 \paragraph{Lower bounding the optimal value.} When running the CoMP algorithm, we at every step $t$ have at our disposal an approximate solution $y^{0,t}$ to the problem of interest (\ref{threepenalties}); $y^{0,t}$ is nothing but the $y^0$-component of the approximate solution $x^t$ generated by CoMP as applied to the saddle point approximation of (\ref{threepenalties}) corresponding to the current value of $\rho_1$,  see (\ref{octeq25}). We have at our disposal also the value $\upsilon(y^{0,t})$ of the objective of (\ref{twopenalties}) at $y^{0,t}$; this quantity is a byproduct of checking whether we should update the current value of $\rho_1$ \footnote{With our implementation, we run this test for both search points and approximate solutions generated by the algorithm}. As a result, we have at our disposal the best found so far value $\upsilon^t=\min_{1\leq\tau\leq t}\upsilon(y^{0,\tau})$, along with the corresponding value $y^{0,t}_*$ of $y^0$: $\upsilon(y^{0,t}_*)=\upsilon^t$. In order to understand how good is the best generated so far approximate solution $y^{0,t}_*$  to the problem of interest, we need to upper bound the quantity $\upsilon^t-\Opt$, or, which is the same, to lower bound $\Opt$. This is a nontrivial task, since the domain of the problem of interest is unbounded, while the usual techniques for online bounding from below the optimal value in a convex minimization problem require the domain to be bounded. We are about to describe a technique for lower bounding $\Opt$ utilizing the structure of (\ref{twopenalties}).
\par
Let $y^0_*$  be an optimal solution to (\ref{twopenalties}) (it clearly exists since $\psi_0\geq0$ and $\lambda,\mu>0$). Assume that at a step $t$ we have at our disposal an upper bound $R=R_t$ on $\|y^0_*\|_1$, and let
$$
R^+=\max[R,\|y^{0,t}\|_1].
$$
Let us look at the saddle point approximation of the problem of interest
\begin{equation}\label{approximation}
\begin{array}{c}
\widehat{\Opt}=\min\limits_{x^1=[y^0;\tau^0;y^1;\tau^1]\in \widehat{X}_1}\max\limits_{x^2\in X_2}\left[\Phi(x^1,x^2):=\psi_0(y^0)+\tau^0+\tau^1+\rho_1\langle y^1-y^0,x^2\rangle\right],\\
X_1=\{[y^0;\tau^0;y^1;\tau^1]:\tau^0\geq{\lambda}\|y^0\|_1,\tau^1\geq\mu\|y^1\|_\nuc\},
\,X_2=\{x^2:\|x^2\|_F\leq1\}.
\end{array}
\end{equation}
associated with current value of $\rho_1$, and let
$$
\bar{X}_1=\{[y^0;\tau^0;y^1;\tau^1]\in X_1:\tau^0\leq {\lambda} R^+,\tau^1\leq \mu R^+\}.
$$
Observe that the point $x^{1,*}=[y^0_*;{\lambda}\|y^0_*\|_1;y^0_*;\mu\|y^0_*\|_\nuc]$ belongs to $\bar{X}_1$ (recall that $\|\cdot\|_\nuc\leq\|\cdot\|_1$) and that
$$
\Opt=\upsilon(y^0_*) \geq \overline{\Phi}(x^{1,*}),\,\,\overline{\Phi}(x^1)=\max_{x^2\in X_2}\Phi(x^1,x^2).
$$
It follows that
$$
\widehat{\Opt}:=\min_{x^1\in \bar{X}_1}\overline{\Phi}(x^1)\leq\Opt.
$$
Further, by Proposition \ref{prop2} as applied to $X^\prime_1=\bar{X}_1$ and $X^\prime_2=X_2$  we have\footnote{note that the latter relation implies that what was denoted by $\widetilde{\Phi}$ in Proposition \ref{prop2} is nothing but $\overline{\Phi}$.}
$$
\overline{\Phi}(x^{1,t})-\widehat{\Opt}\leq \Res(\bar{X}_1\times X_2\big|\cI_t,\lambda^t),
$$
where $\cI_t$ is the execution protocol generated by CoMP {\sl as applied to the saddle point problem} (\ref{approximation}) (i.e., since the last restart preceding step $t$ till this step), and $\lambda^t$ is the associated accuracy certificate. We conclude that
$$
\ell_t:=\overline{\Phi}(x^{1,t})-\Res(\bar{X}_1\times X_2\big|\cI_t,\lambda^t)\leq\widehat{\Opt}\leq\Opt,
$$
and $\ell_t$ is easy to compute (since the resolution is just the maximum of a readily given by $\cI_t,\lambda^t$ affine function over $\bar{X}_1\times X_2$). Setting $\upsilon_t=\max_{\tau\leq t} \ell_\tau$, we get nondecreasing with $t$  lower bounds on $\Opt$. Note that this component of our lower bounding is independent of the particular structure of $\psi_0$.\par
It remains to explain how to get an upper bound $R$ on $\|y^0_*\|_1$, and this is where the special structure of $\psi_0(y)={1\over 2}\|P_\Omega y-b\|_2^2$ is used. Recalling that $b\in\bR^M$, let us set
$$
\vartheta(r) =\min_{v\in\bR^M}\left\{{1\over 2}\|v-b\|_2^2:\|v\|_1\leq r\right\},\,\,r\geq 0,
$$
It is immediately seen that replacing the entries in $b$ by their magnitudes, $\vartheta(\cdot)$ remains intact, and that for $b\geq0$ we have
$$
\vartheta(r)=\min_{v\in\bR^M}\left\{{1\over 2}\|v-b\|_2^2: v\geq0,\sum_iv_i\leq r\right\},
$$
so that $\vartheta(\cdot)$ is an easy to compute nonnegative and nonincreasing convex function of $r\geq0$. Now, by definition of $P_\Omega$, the function $\vartheta^+(\|y^0\|_1)$ where
$$
\vartheta^+(r)= {\lambda} r + \vartheta(r)
$$
 is a lower bound on $\upsilon(y^0)$. As a result, given an upper bound $\upsilon^t$ on $\Opt=\upsilon(y_*)$, the easy-to-compute quantity
$$
R_t:=\max\{r: \vartheta^+(r)\leq \upsilon^t\}
$$
is an upper bound on $\|y^0_*\|_1$. Since $\upsilon^t$ is nonincreasing in $t$, $R_t$ is nonincreasing in $t$ as well.

\paragraph{Generating the data.} In the experiments to be reported, the data of (\ref{twopenalties}) were generated as follows. Given $n$, we build ``true'' $n\times n$ matrix
 $y_{\#}=\sum_{i=1}^ke_if_i^T$, with $k=\lfloor n/4\rfloor$ and vectors $e_i,f_i\in\bR^n$  sampled, independently of each other, as follows: we draw a vector from the standard Gaussian distribution $\cN(0,I_n)$, and then zero out part of the entries, with probability of replacing a particular entry with zero selected in such a way that the sparsity of $y_{\#}$ is about a desired level (in our experiments, we wanted $y_{\#}$ to have about 10\% of nonzero entries).  The set $\Omega$ of ``observed cells'' was built at random, with probability 0.25 for a particular cell to be in $\Omega$. Finally, $b$ was generated
 as $P_\Omega(y_{\#}+\sigma \xi)$, where the entries of $\xi\in\bR^{n\times n}$ were independently of each other drawn from the  standard Gaussian distribution, and
 $$
 \sigma =0.1{\sum_{i,j}|[y_{\#}]_{ij}|\over n^2}.
 $$
 We used ${\lambda}=\mu=10\sigma$.\footnote{If the goal of solving (\ref{twopenalties}) were to recover $y_{\#}$, our ${\lambda}$ and $\mu$ would, perhaps, be too large. Our goal, however, was solving (\ref{twopenalties}) as an ``optimization beast,'' and we were interested in ``meaningful'' contribution of $\Psi_0$ and $\Psi_1$ to the objective of the problem, and thus in not too small ${\lambda}$ and $\mu$.}
 Finally, our guess for the Frobenius norm of the optimal solution to (\ref{twopenalties}) is defined as follows. Note that the quantity $\|b\|_2^2-M\sigma^2$ is an estimate of $\|P_\Omega y_{\#}\|_2^2$. We define the estimate $D$ of $D_*:=\|y_*\|_F$
``as if'' the optimal solution were $y_{\#}$, and all entries of $y_{\#}$ were of the same order of magnitude
\[
 D=\sqrt{{n^2\over M}\max[\|b\|_2^2-M\sigma^2,1]},\,\,M=\Card(\Omega).
\]

 \paragraph{Numerical results.} The results of the first series of experiments are presented in Table \ref{tab01}. The comments are as follows.
\par
In the ``small'' experiment $(n=128$, the largest $n$ where we were able to solve (\ref{twopenalties}) in a reasonable time by {\tt CVX} \cite{CVX} {using the state-of-the-art   {\tt mosek} \cite{Mosek} Interior-Point solver} and thus knew the ``exact'' optimal value), CoMP exhibited fast convergence: relative accuracies 1.1e-3 and 6.2e-6 are achieved in 64 and 4096 steps (1.2 sec and 74.9 sec,
     respectively, as compared to 4756.7 sec taken by {\tt CVX}).
\par In larger experiments ($n=512$ and $n=1024$, meaning design dimensions  262,144 and  1,048,576, respectively), the running times look moderate, and the convergence pattern of the CoMP still looks promising\footnote{Recall that we do not expect linear convergence, just $O(1/t)$ one.}. Note that our lower bounding, while somehow working, is very conservative: it overestimates the ``optimality gap'' $\upsilon^t-\upsilon_t$ by 2-3 orders of magnitude  for moderate and large values of $t$ in the $128\times 128$ experiment. More accurate performance evaluation would require a less conservative lower bounding of the optimal value  (as of now, we are not aware of any alternative).
\par
In the second series of experiments, the data of (\ref{twopenalties}) were generated in such a way that the true optimal solution and optimal value to the problem were known from the very beginning. To this end we take as $\Omega$ the collection of all cells of an $n\times n$ matrix, which, via optimality conditions, allows to select $b$ making our ``true'' matrix $y_{\#}$ the optimal solution to (\ref{twopenalties}). The results are presented in Table \ref{tab02}.
\par
It should be mentioned that in these experiments the value of $\rho_1$ resulting in negligibly small, as compared to $\epsilon_2$, values of $\epsilon_1$ in (\ref{epsilonone}) was found in the first 10-30 steps of the algorithm, with no restarts afterwards.

\begin{table}
\begin{center}
{\scriptsize
\begin{tabular}{|c||r|r|r|r|r|r|r|r|r|r|}
\hline
$t$&8&16&32&64&128&256&512&1024&2048&4096\\
\hline\hline
CPU, sec&0.1 &   0.3 &   0.6 &   1.2 &   2.3 &   4.7 &   9.4 &  18.7 &  37.5  & 74.9\\
\hline\hline
$\upsilon^t-\Opt$&2.0e-2 &1.8e-2& 1.8e-2 &1.4e-2 &5.3e-3& 5.0e-3 &1.3e-3& 7.8e-4& 3.2e-4 &8.3e-5\\
\hline
$\upsilon^t-\upsilon_t$&4.8e0& 4.5e0 &4.2e0& 3.7e0& 2.1e0& 6.3e-1& 2.1e-1& 1.3e-1& 6.0e-2& 3.4e-2\\
\hline\hline
${\upsilon^t-\Opt\over\Opt}$&1.5e-3& 1.3e-3& 1.3e-3& 1.1e-3& 4.0e-4& 3.7e-4& 9.5e-5& 5.8e-5& 2.4e-5& 6.2e-6\\
\hline
${\upsilon^t-\upsilon_t\over\upsilon_{4096}}$&3.6e-1& 3.4e-1& 3.2e-1& 2.8e-1& 1.5e-1& 4.7e-2& 1.6e-2& 9.4e-3& 4.5e-3& 2.6e-3\\
\hline\hline
${\upsilon^1-\Opt\over\upsilon^t-\Opt}$&4.8e1& 5.4e1& 5.4e1& 6.7e1& 1.8e2& 1.9e2& 7.5e2 &1.2e3& 2.9e3& 1.1e4\\
\hline
${\upsilon^1-\upsilon_1\over \upsilon^t-\upsilon_t}$&3.0e0& 3.2e0 &3.7e0 &3.9e0 &6.9e0 &2.3e1& 6.7e1& 1.1e2& 2.4e2& 4.1e2\\
\hline\hline
\multicolumn{10}{c}{}\\
\multicolumn{10}{c}{(a) $n=128$, $\Opt=13.28797$ (CVX CPU 4756.7 sec)}\\
\multicolumn{10}{c}{}\\
\cline{1-10}
$t$&8& 16& 32& 64& 128 &256& 512& 1024&2048&\multicolumn{1}{|c}{}\\
\cline{1-10}
CPU, sec&3.7&    7.5 &  15.0 &  29.9 &  59.8 & 119.6 & 239.2 & 478.4&992.0&\multicolumn{1}{|c}{}\\
\cline{1-10}
$\upsilon^t-\upsilon_t$&4.4e1& 4.4e1& 4.3e1& 4.2e1& 4.1e1& 3.7e1& 2.3e1& 1.2e1& 5.1e0&\multicolumn{1}{|c}{}\\
\cline{1-10}
${\upsilon^t-\upsilon_t\over\upsilon_{1024}}$&2.4e-1& 2.4e-1& 2.4e-1& 2.4e-1& 2.2e-1& 2.0e-1& 1.3-1& 6.4e-2& 2.8e-2&\multicolumn{1}{|c}{}\\
\cline{1-10}
${\upsilon^1-\upsilon_1\over\upsilon^t-\upsilon_t}$&4.4e0& 4.4e0& 4.5e0& 4.6e0& 4.8e0& 5.5e0& 8.5e0& 1.7e1&3.8e1&\multicolumn{1}{|c}{}\\
\cline{1-10}
\multicolumn{11}{c}{}\\
\multicolumn{10}{c}{(b) $n=512$, $\upsilon_{2048}=175.445 \leq \Opt\leq \upsilon^{2048}=180.503$ (CVX not tested)}&\multicolumn{1}{c}{}\\
\multicolumn{11}{c}{}\\
\cline{1-9}
$t$&8& 16& 32& 64& 128 &256& 512& 1024&\multicolumn{2}{|c}{}\\
\cline{1-9}
CPU, sec& 23.5&   46.9 &  93.8 & 187.6 & 375.3 & 750.6& 1501.2& 3002.3&\multicolumn{2}{|c}{}\\
\cline{1-9}
$\upsilon^t-\upsilon_t$&1.5e2 &1.5e2&1.3e2& 1.2e2& 1.1e2& 8.0e1 &1.6e1& 5.4e0&\multicolumn{2}{|c}{}\\
\cline{1-9}
${\upsilon^t-\upsilon_t\over\upsilon_{1024}}$&2.4e-1& 2.2e-1& 2.2e-1& 1.9e-1& 1.7e-01& 1.2e-1& 2.4e-2& 8.1e-3&\multicolumn{2}{|c}{}\\
\cline{1-9}
${\upsilon^1-\upsilon_1\over \upsilon^t-\upsilon_t}$&4.6e0 &4.8e0& 5.3e0 &5.7e0& 6.3e0& 8.9e0 &4.5e1 &1.3e2&\multicolumn{2}{|c}{}\\
\cline{1-9}
\multicolumn{9}{c}{}\\
\multicolumn{9}{c}{(c) $n=1024$, $\upsilon_{1024}=655.422\leq\Opt\leq \upsilon^{1024}=660.786$ (CVX not tested)}\\
\end{tabular}}
\end{center}
\caption{\label{tab01} Composite Mirror Prox algorithm on problem (\ref{twopenalties}) with $n\times n$ matrices.  $\upsilon^t$ are the best values of $\upsilon(\cdot)$, and $\upsilon_t$ are lower bounds on the optimal value found in course of $t$ steps.
Platform:  {\tt MATLAB} on
3.40 GHz Intel Core i7-3770 desktop with 16 GB RAM, 64 bit Windows 7.}
\end{table}

 \paragraph{Remarks.} For the sake of simplicity, so far we were considering problem (\ref{twopenalties}), where minimization is  carried out over $y^0$ running through the entire space $\bR^{n\times n}$ of $n\times n$ matrices. What happens if we restrict $y^0$ to reside in a given closed convex domain $Y_0$?\par
 It is immediately seen that the construction we have presented can be straightforwardly modified for the cases when $Y_0$ is a centered at the origin ball in the Frobenius or $\|\cdot\|_1$ norm, or the intersection of such a set with the space of symmetric $n\times n$ matrices. We could also handle the case when $Y_0$ is the centered at the origin nuclear norm ball (or intersection of this ball with the space of symmetric matrices, or with the cone of positive semidefinite symmetric matrices), but to this end one needs to ``swap the penalties'' -- to write the representation (\ref{octeq23c}) of problem (\ref{twopenalties}) as
 $$
 \begin{array}{c}
 \min\limits_{{\{y^k;\tau^k]\}_{k=0}^1\atop\in Y_0^+\times Y_1^+}}\bigg\{\Upsilon(y^0,y^1,\tau^0,\tau^1):=\underbrace{{1\over 2}\|P_\Omega y^0-b\|_2^2}_{\psi_0(y^0)} +\tau^0+\tau^1:y^0=y^1\bigg\},\\
 Y_0^+=\{[y^0;\tau^0]: y^0\in Y_0,\tau^0\geq\mu\|y^0\|_\nuc\},\,\,Y_1^+=\{[y^1;\tau^1]: y^1\in Y_1,\tau^1\geq{\lambda}\|y^1\|_1\},\\
  \end{array}
 $$
 where $Y_1\supset Y_0$ ``fits'' $\|\cdot\|_1$ (meaning that we can point out a d.g.f. $\omega_1(\cdot)$ for $Y_1$ which, taken along with $\Psi_1(y^1)={\lambda}\|y^1\|_1$, results in easy-to-solve auxiliary problems (\ref{octeq20})). We can take, e.g. $\omega_1(y^1)={1\over 2}\|y^1\|_F^2$ and define $Y_1$ as the entire space, or a centered at the origin  Frobenius/$\|\cdot\|_1$ norm ball large enough to contain $Y_0$.

\begin{table}
\begin{center}
{\scriptsize
\begin{tabular}{|c||r|r|r|r|r|r|r|r|}
\hline
$t$&1&7&8&12&128&256&512&1024\\
\hline
CPU, sec&1.3&8.3 &9.3&11.0&65.9&125.0&244.7&486.0\\
\hline
$\upsilon^t-\Opt$&92.9&1.58&0.30&0.110&0.095&0.076&0.069&0.069\\
\hline
$\upsilon^t-\upsilon_t$&700.9&92.4&69.5&54.6&52.8&44.2&21.2&3.07\\
\hline
${\upsilon^t-\Opt\over\Opt}$&0.153&2.6e-3&5.0e-4&1.8e-4&1.6e-4&1.3e-4&1.1e-4&1.1e-4\\
\hline
${\upsilon^t-\upsilon_t\over\Opt}$&1.153&0.152&0.114&0.090&0.087&0.073&0.035&0.005\\
\hline
\multicolumn{9}{c}{~}\\
\multicolumn{9}{c}{(a) $n=512$, $\Opt=607.9854$}\\
\multicolumn{9}{c}{}\\
\cline{1-7}
$t$&1&7&8&128&256&512&\multicolumn{2}{c}{}\\
\cline{1-7}
CPU, sec&8.9 &48.1&51.9&392.7&752.1&1464.9&\multicolumn{2}{c}{}\\
\cline{1-7}
$\upsilon^t-\Opt$&371.4&3.48&0.21&0.21&0.19&0.16\\
\cline{1-7}
$\upsilon^t-\upsilon_t$&2772&241.7&201.2&147.3&146.5&122.9\\
\cline{1-7}
${\upsilon^t-\Opt\over\Opt}$&0.154&1.5e-3&9e-5&9e-5&8e-5&7e-5\\
\cline{1-7}
${\upsilon^t-\upsilon_t\over\Opt}$&1.155&0.101&0.084&0.061&0.061&0.051\\
\cline{1-7}
\multicolumn{7}{c}{}\\
\multicolumn{7}{c}{(b) $n=1024$, $\Opt=2401.168$}&\multicolumn{2}{c}{}\\
\multicolumn{7}{c}{}\\
\end{tabular}}
\end{center}
\caption{\label{tab02} Composite Mirror Prox algorithm on problem (\ref{twopenalties}) with $n\times n$ matrices and known optimal value $\Opt$. $\upsilon^t$ are the best values of $\upsilon(\cdot)$, and $\upsilon_t$ are lower bounds on the optimal value found in course of $t$ steps.
Platform:  {\tt MATLAB} on 3.40 GHz Intel Core i7-3770 desktop with 16 GB RAM, 64 bit Windows 7.}
\end{table}

\subsubsection{Image decomposition}

\paragraph{Problem of interest.} In the experiments to be reported, we applied the just outlined approach
to {Example 2, that is, to the problem}
\be
\Opt=\min\limits_{y^1,y^{2},y^{3}\in\bR^{n\times n}}\left\{\|A (y^{1}+y^{2}+y^{3}) -b\|_2+\mu_1\|y^{1}\|_\nuc +\mu_2\|y^{2}\|_1+\mu_3 \|y^{3}\|_\TV\right\}.
\ee{example2}
{where $A(y):\,\bR^{n\times n}\to \bR^{M}$ is a given linear mapping.}

\paragraph{Problem reformulation.} We first rewrite \rf{example2} as
a saddle point optimization  problem
\be
\Opt=\min\limits_{y^{1},y^{2},y^{3}\in\bR^{n\times n}}\left\{\|A(y^{1}+y^{2}+y^{3}) -b\|_2+\mu_1\|y^{1}\|_\nuc +\mu_2 \|y^{2}\|_1+\mu_3\|T y^{3}\|_1\right\}\nn
=\min\limits_{y^{1},y^{2},y^{3}\in\bR^{n\times n}}\left\{\max_{\|z\|_2\le 1}\langle z,A(y^{1}+y^{2}+y^{3}) -b\rangle
+\mu_1\|y^{1}\|_\nuc +\mu_2 \|y^{2}\|_1+\mu_3\|T y^{3}\|_1\right\},
\ee{threepenalties}
where $T:\;\bR^{n\times n}\to \bR^{2n(n-1)}$ is the mapping  $y\mapsto T y=\left[
\begin{array}{c}\{(\nabla_i y)_{n(j-1)+i}\}_{i=1,...,n-1,\,j=1,...,n}\\
\{(\nabla_j y)_{n(i-1)+j})\}_{i=1,...,n,\,j=1,...,n-1}\end{array}\right]$.\\

Next we rewrite \rf{threepenalties} as a linearly constrained saddle-point problem with ``simple'' penalties:
\bse
\Opt&=&\min\limits_{
{y^3\in Y_3\atop [y^k;\tau_k]\in Y^+_k,\, 0\leq k\leq 2} }\max_{z\in Z}\left\{\langle z,A(y^1+y^2+y^3) -b\rangle
+\tau_1 +\tau_2+\tau_0,\;y^0=Ty^{3}\right\},
\ese
where \bse
Y_0^+&=&\{[y^0;\tau_0]:\, y^0\in Y_0=\bR^{2n(n-1)}: \|y^{0}\|_{1}\leq \tau_0/\mu_3\},\;\;\\
Y_1^+&=&\{[y^1;\tau_1]:\, y^1\in Y_1=\bR^{n\times n}: \|y^{1}\|_{\nuc}\leq \tau_1/\mu_1\},\;\;\\Y_2^+&=&\{[y^2;\tau_2]:\, y^2\in Y_2=\bR^{n\times n}: \|y^{2}\|_{1}\leq \tau_2/\mu_2\}\\
Y_3&=&\bR^{n\times n},
\;\;Z=\{z\in \bR^{M}: \|z\|_2\leq 1\},
\ese
and further approximate the resulting problem with its penalized version:
\be
\widehat{\Opt}&=\min\limits_{
{y^3\in Y_3\atop [y^k;\tau_k]\in Y^+_k,\, 0\leq k\leq 2} }\max\limits_{z\in Z\atop w\in W }
&\left\{\begin{array}{l}\langle z,A(y^1+y^2+y^3) -b\rangle
\\
+\tau_1 +\tau_2+\tau_0+\rho \langle w,y^0-Ty^{3}\rangle\end{array}\right\},
\ee{3pensaddle}
with
$$
W=\{w\in \bR^{2n(n-1)},\,\|w\|_2\leq 1\}.
$$
Note that the function $\psi(y^1,y^2,y^3):=\|A(y^1+y^2+y^3)-b\|_2=\max_{\|z\|_2\leq 1} \langle z,\,A(y^1+y^2+y^3)-b\rangle $ is Lipschitz continuous in $y^3$ with respect to the Euclidean norm on $\bR^{n\times n}$ with corresponding Lipschitz constant $G=\|A\|_{2,2}$, which is the spectral norm (the principal singular value) of $A$. Further, $\Psi(y^0)=\mu_3\|y^0\|_1$ is Lipschitz-continuous in $y^0$ with respect to the Euclidean norm on $\bR^{2n(n-1)}$ with the Lipschitz constant $H\leq \mu_3\sqrt{2n(n-1)}$. With the help of the result of Proposition \ref{pro:correction} we conclude that to ensure the ``exact penalty'' property it suffices to choose $\rho\geq \|A\|_{2,2}+\mu_3\sqrt{2n(n-1)}$.
Let us denote
\[
U=\left\{\begin{array}{c}u=[y^0;...; y^3; z; w]: \;y^k\in Y^k,\,0\leq k\leq 3,\\
z\in \bR^M,\,\|z\|_2\leq 1,\,w\in \bR^{2n(n-1)},\,\|w\|_2\leq 1\end{array} \right\}
\]
We equip the embedding space $E_u$ of $U$ with the norm
\[
\|u\|=\left(\alpha_0\|y^0\|_2^2+\sum_{k=1}^3 \alpha_k\|y^k\|^2_2+\beta\|z\|_2^2+\gamma\|w\|^2_2
\right)^{1/2},
\]
and $U$ with the proximal setup $(\|\cdot\|,\,\omega(\cdot))$ with
\[
\omega(u)=
{\alpha_0\over 2}\|y^0\|_2^2+\sum_{k=1}^3 {\alpha_k\over 2}\|y^k\|^2_2+\half \|z\|_2^2+\half \|w\|_2^2
\]

\paragraph{Implementing the CoMP algorithm.} When implementing the CoMP algorithm, we use the above proximal setup with adaptive aggregation parameters $\alpha_0=\cdots=\alpha_4= 1/D^2$ where $D$ is our guess for the upper bound of $||y_*||_2$, that is, whenever the norm of the current solution exceeds $20\%$ of the guess value, we increase $D$ by factor $2$ and update the scales accordingly.  The penalty $\rho$ and stepsizes $\gamma_t$ are adjusted dynamically the same way as explained in the last experiment.

\paragraph{Numerical results.} In the first series of experiments, we build  the $n\times n$ observation matrix $b$ by first generating a random matrix with rank $r=\lfloor\sqrt{n}\rfloor$ and another random matrix with sparsity $p=0.01$, so that the  observation matrix is a sum of these two matrices and of random noise of level $\sigma=0.01$; we take $y\mapsto Ay$ as the identity mapping. We use $\mu_1=10\sigma, \mu_2=\sigma,\mu_3=\sigma$.
The very preliminary results of this series of experiments are presented in Table \ref{tab03}.
Note that unlike the matrix completion problem, discussed in section \ref{sec:mcomplet}, here we are not able to generate the problem with known optimal solutions.  Better performance evaluation would require good lower bounding of the true optimal value, which is however problematic due to unbounded problem domain.

\begin{table}[f]
\begin{center}
{\scriptsize
\begin{tabular}{|c||r|r|r|r|r|r|r|r|r|}
\hline				
$t$	&8	&16	&32	&64	&128	&256	&512	&1024	&2048\\
\hline\hline
CPU, sec	&0.1	&0.2	&0.4	&0.8	&1.6	&3.1	&6.3	&12.6	&25.2\\
\hline\hline
$\upsilon_t-\upsilon_{2048}$	&1.5e1	&2.8e0	&6.2e-1	&2.3e-1	&1.1e-1	&4.2e-2	&1.5e-2	&4.4e-3	&0.0e0\\
\hline
$\frac{\upsilon_t-\upsilon_{2048}}{\upsilon_{2048}}$	&9.5e-1	&1.8e-1	&4.0e-2	&1.5e-2	&7.0e-3	&2.7e-3	&9.9e-4	&2.8e-4	&0.0e0\\
\hline\hline
$\upsilon_t-\Opt$	&1.5e1	&2.8e0	&6.2e-1	&2.3e-1	&1.1e-1	&4.5e-2	&1.8e-2	&6.6e-3	&2.2e-3\\
\hline
$\frac{\upsilon_t-\Opt}{\Opt}$	&9.5e-1	&1.8e-1	&4.0e-2	&1.5e-2	&7.1e-3	&2.9e-3	&1.1e-3	&4.2e-4	&1.4e-4\\
\hline
\multicolumn{10}{c}{}\\
\multicolumn{10}{c}{(a) $n=64$, $\Opt=15.543$ (CVX CPU 4525.5 sec)}\\
\multicolumn{10}{c}{}\\
\hline
$t$	&8	&16	&32	&64	&128	&256	&512	&1024	&2048\\
\hline\hline
CPU, sec	&6.2	&12.3	&24.7	&49.3	&98.6	&197.2	&394.4	&788.9	&1577.8\\
\hline\hline
$\upsilon_t-\upsilon_{2048}$	&1.1e2	&5.8e1	&2.7e1	&1.3e1	&6.2e0	&2.9e0	&1.2e0	&3.9e-1	&0.0e0\\
\hline
$\frac{\upsilon_t-\upsilon_{2048}}{\upsilon_{2048}}$	&9.0e-1	&4.9e-1	&2.3e-1	&1.1e-1	&5.2e-2	&2.5e-2	&1.0e-2	&3.3e-3	&0.0e0\\
\hline
\multicolumn{10}{c}{}\\
\multicolumn{9}{c}{(b) $n=512$  (CVX not tested)}\\
\end{tabular}}
\end{center}
\caption{\label{tab03} Composite Mirror Prox algorithm on problem (\ref{example2}) with $n\times n$ matrices.  $\upsilon^t$ are the best values of $\upsilon(\cdot)$ in course of $t$ steps.
Platform: {\tt MATLAB} on Intel i5-2400S @2.5GHz CPU, 4GB RAM, 64-bit Windows 7.}
\end{table}

In the second experiment we implemented the CoMP algorithm to decompose real images and extract the underlying low rank{/sparse singular distortion/smooth background} {components}. The purpose of these experiments is to illustrate how the algorithm performs with the choice of small regularization parameters which is meaningful from the point of view of applications to image recovery. Image decomposition results for two images are provided on figures \ref{checker} and \ref{building}.
On figure \ref{checker} we present the decomposition of the observed image of size $256\times 256$. We apply the model (\ref{threepenalties}) with regularization parameters $\mu_1=0.03,\mu_2=0.001,\mu_3=0.005$. We run  $2\,000$ iterations  of CoMP (total of $393.5$ sec {\tt MATLAB}, Intel i5-2400S @2.5GHz CPU). The first component $y_1$ has approximate rank $\approx1$; the relative error of the reconstruction  $\|y_1+y_2+y_3-b\|_2/\|b\|_2\approx 2.8\times 10^{-4}$.

Figure \ref{building} shows the decomposition of the observed image of size $480\times 640$ after $1\,000$ iterations of CoMP ({CPU $873.6$ sec}). The regularization parameters of the model  (\ref{threepenalties}) were set to $\mu_1=0.06,\mu_2=0.002,\mu_3=0.005$.
The  relative error of the reconstruction  $\|y_1+y_2+y_3-b\|_2/\|b\|_2\approx 8.4\times 10^{-3}$.\\

\begin{figure}[f]
\centering
\subfloat[observation $b$]{
                \includegraphics[trim=1cm 1cm 1cm 4mm, clip, width=0.3\textwidth]{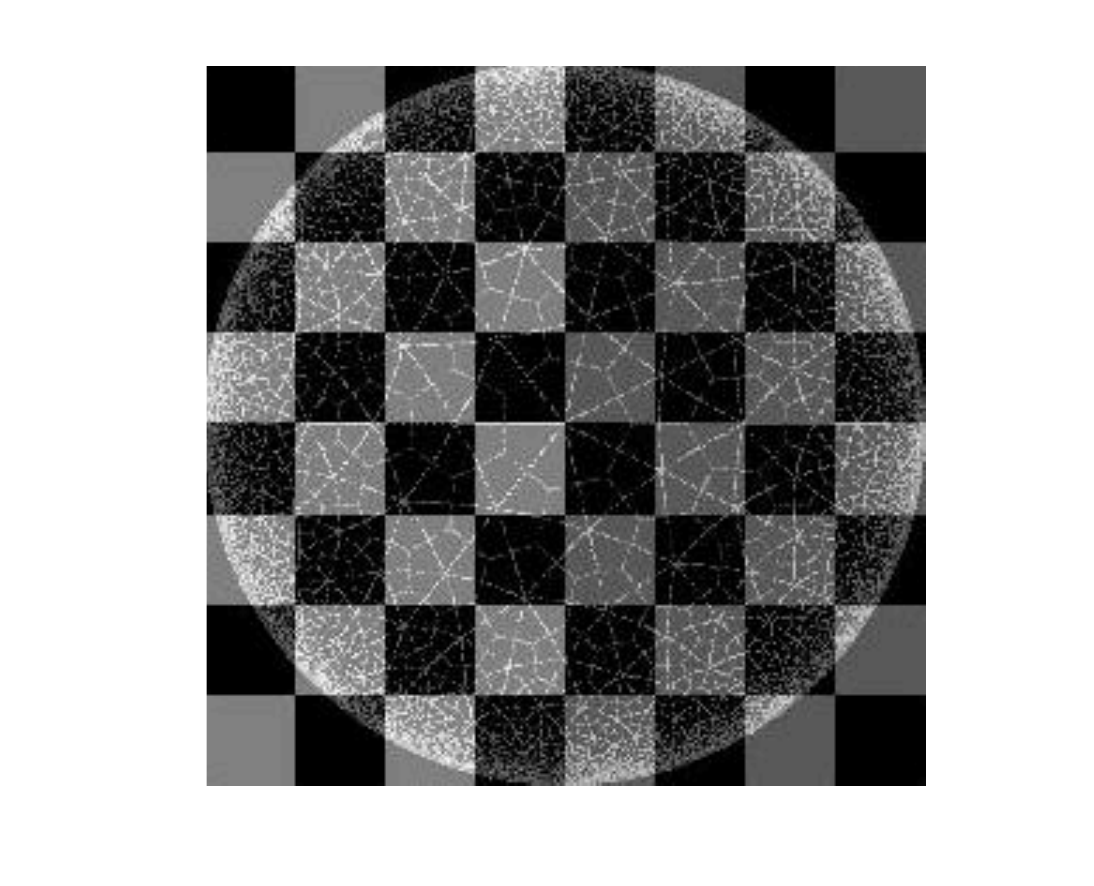}
                \label{fig1:obs}}
\subfloat[recovery $y_1+y_2+y_3$]{
                \includegraphics[trim=1cm 1cm 1cm 4mm, clip, width=0.3\textwidth]{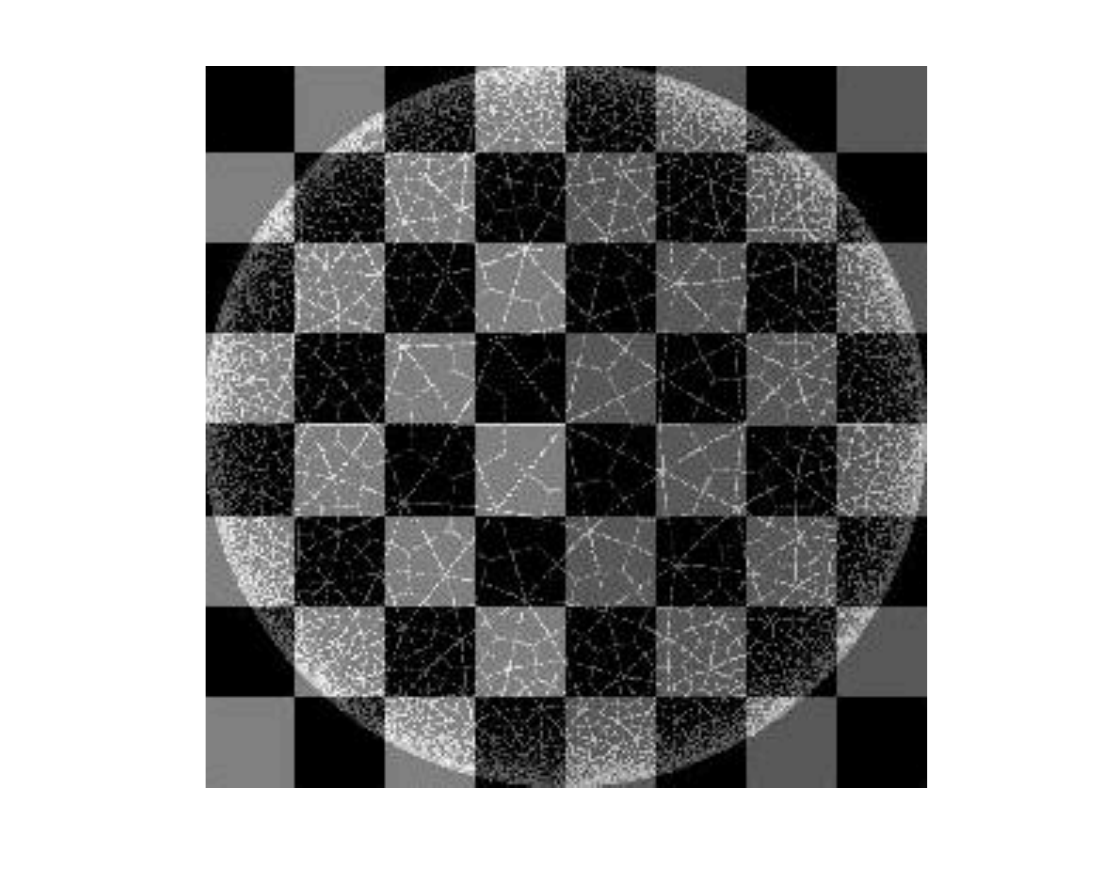}
                \label{fig1:recovery}}

\subfloat[low-rank component]{
         \includegraphics[trim=1cm 1cm 1cm 4mm, clip,width=0.3\textwidth]{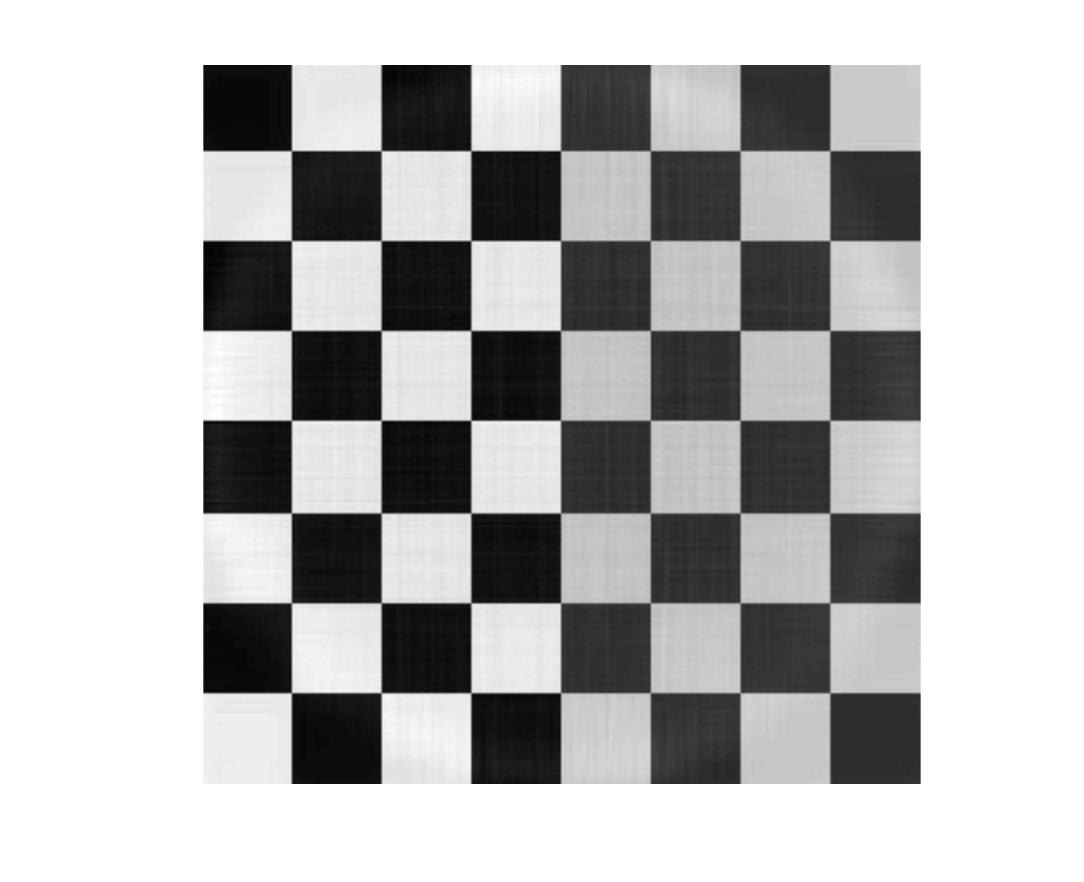}
                \label{fig1:lowrank}}
\subfloat[sparse component]{
           \includegraphics[trim=1cm 1cm 1cm 4mm, clip,width=0.3\textwidth]{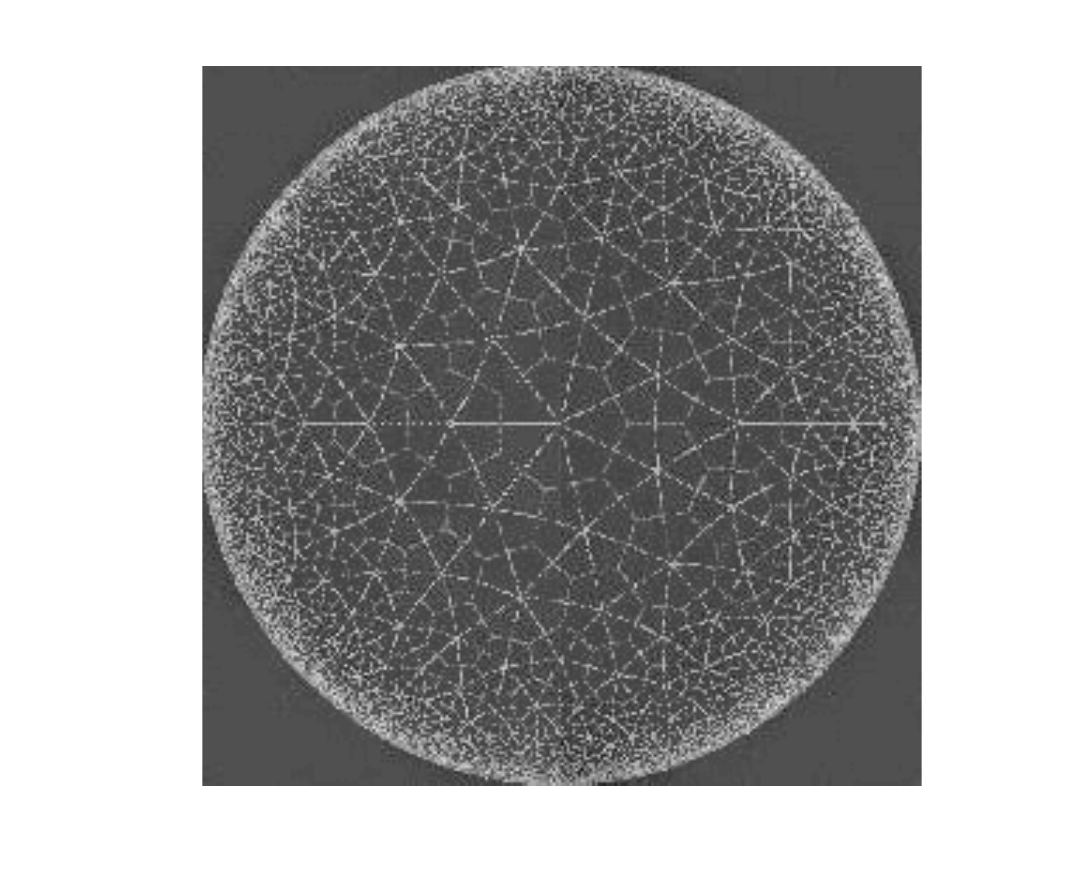}
                \label{fig1:sparse}}
\subfloat[smooth component]{
        \includegraphics[trim=1cm 1cm 1cm 4mm, clip,width=0.3\textwidth]{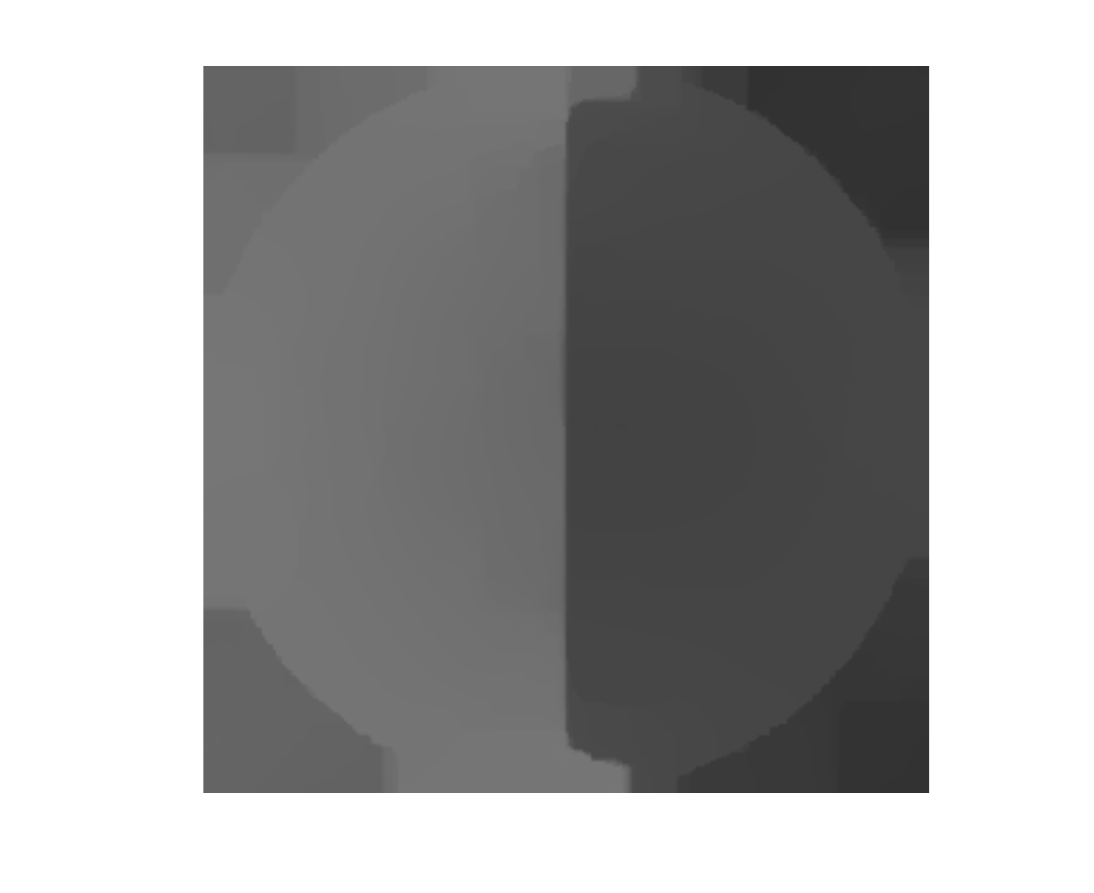}
                \label{fig1:smooth}}
\caption{Observed and reconstructed images (size $256\times256$).
}
\label{checker}
\end{figure}
\begin{figure}[f]
\centering
\subfloat[observation $b$]{
                \includegraphics[trim=1cm 1.2cm 1cm 8mm, clip, width=0.45\textwidth]{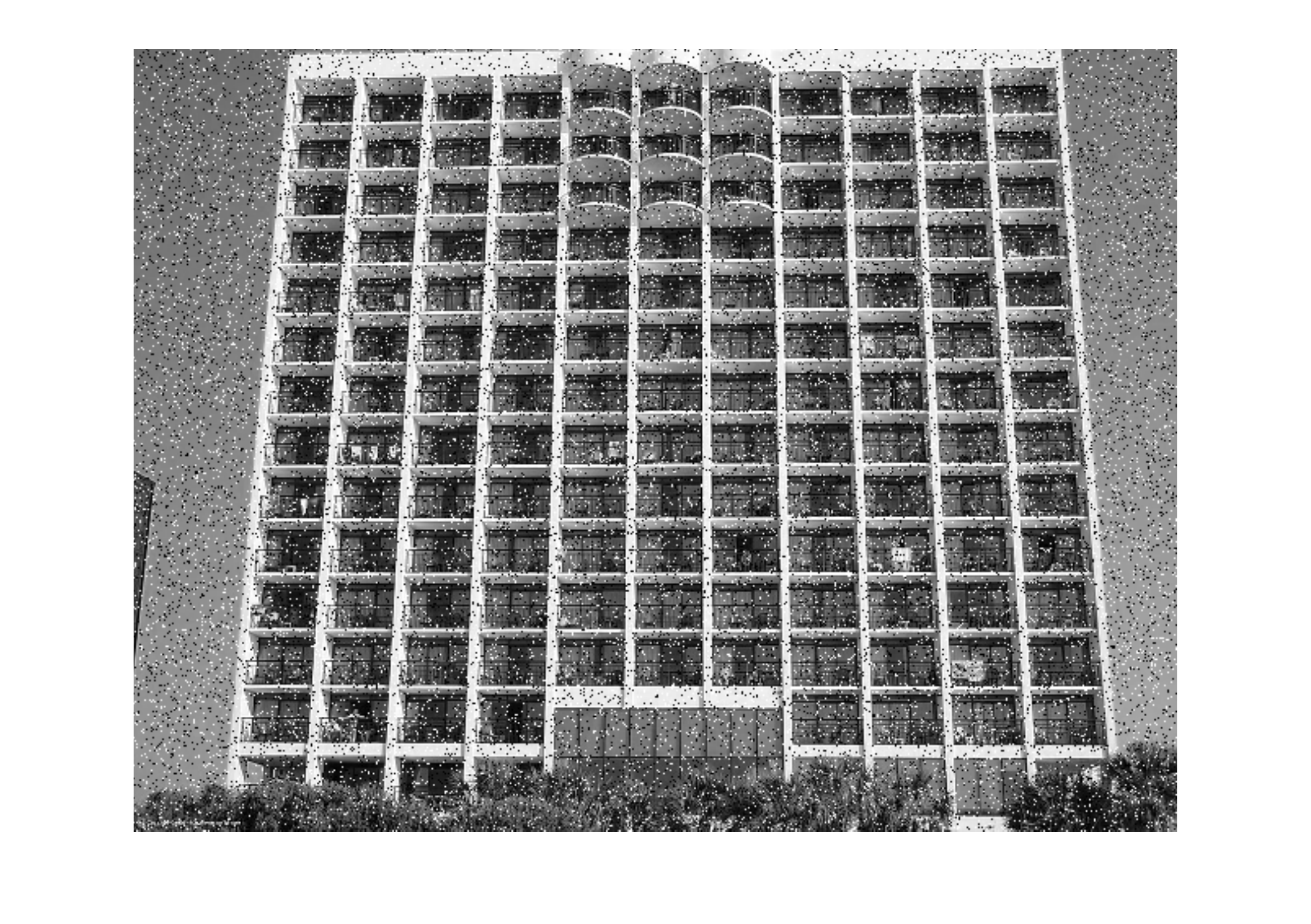}
                \label{fig2:obs}}
\subfloat[low-rank component]{
         \includegraphics[trim=1cm 1.2cm 1cm 8mm, clip,width=0.45\textwidth]{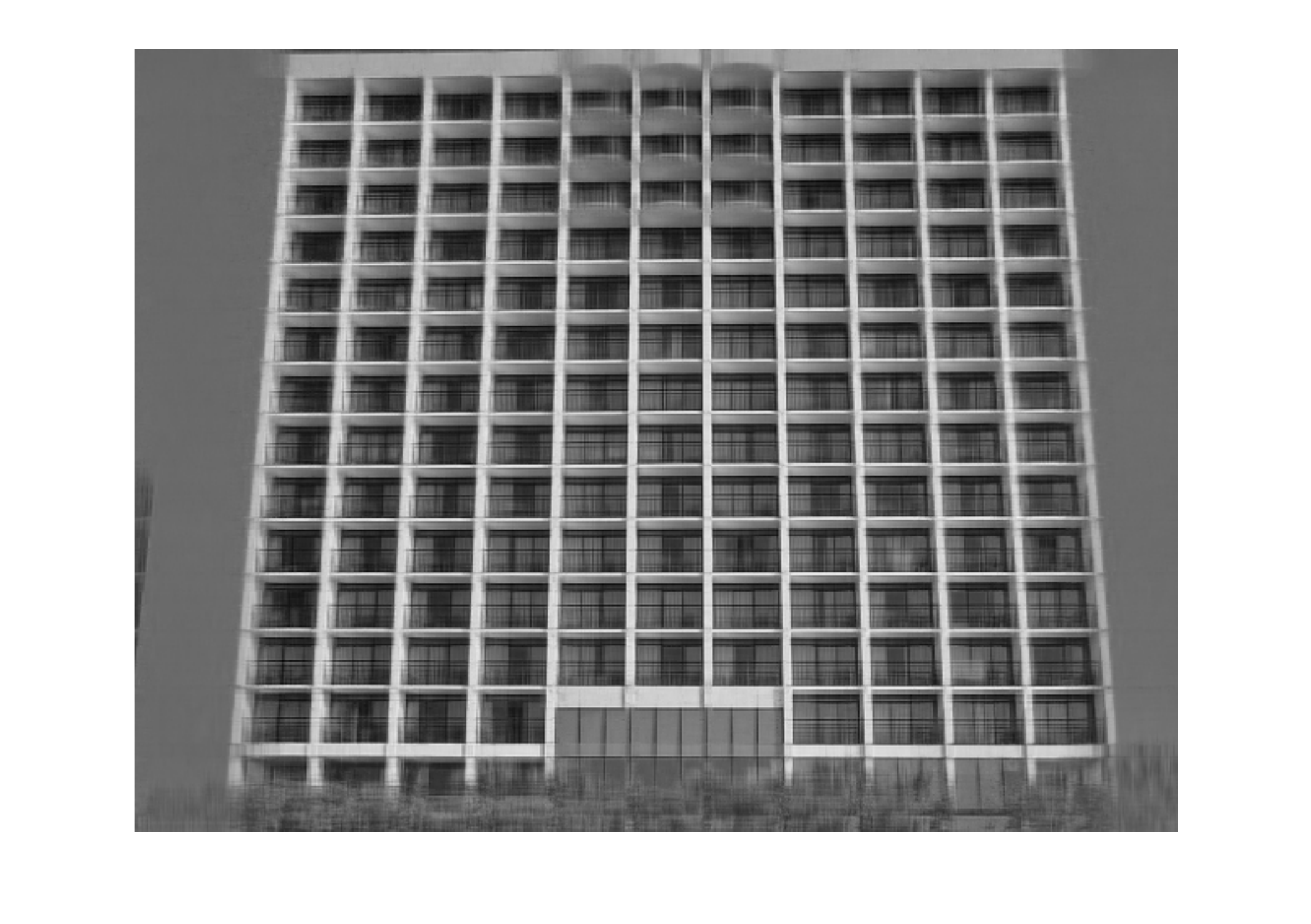}
                \label{fig2:lowrank}}

\subfloat[sparse component]{
           \includegraphics[trim=1cm 1.2cm 1cm 8mm, clip,width=0.45\textwidth]{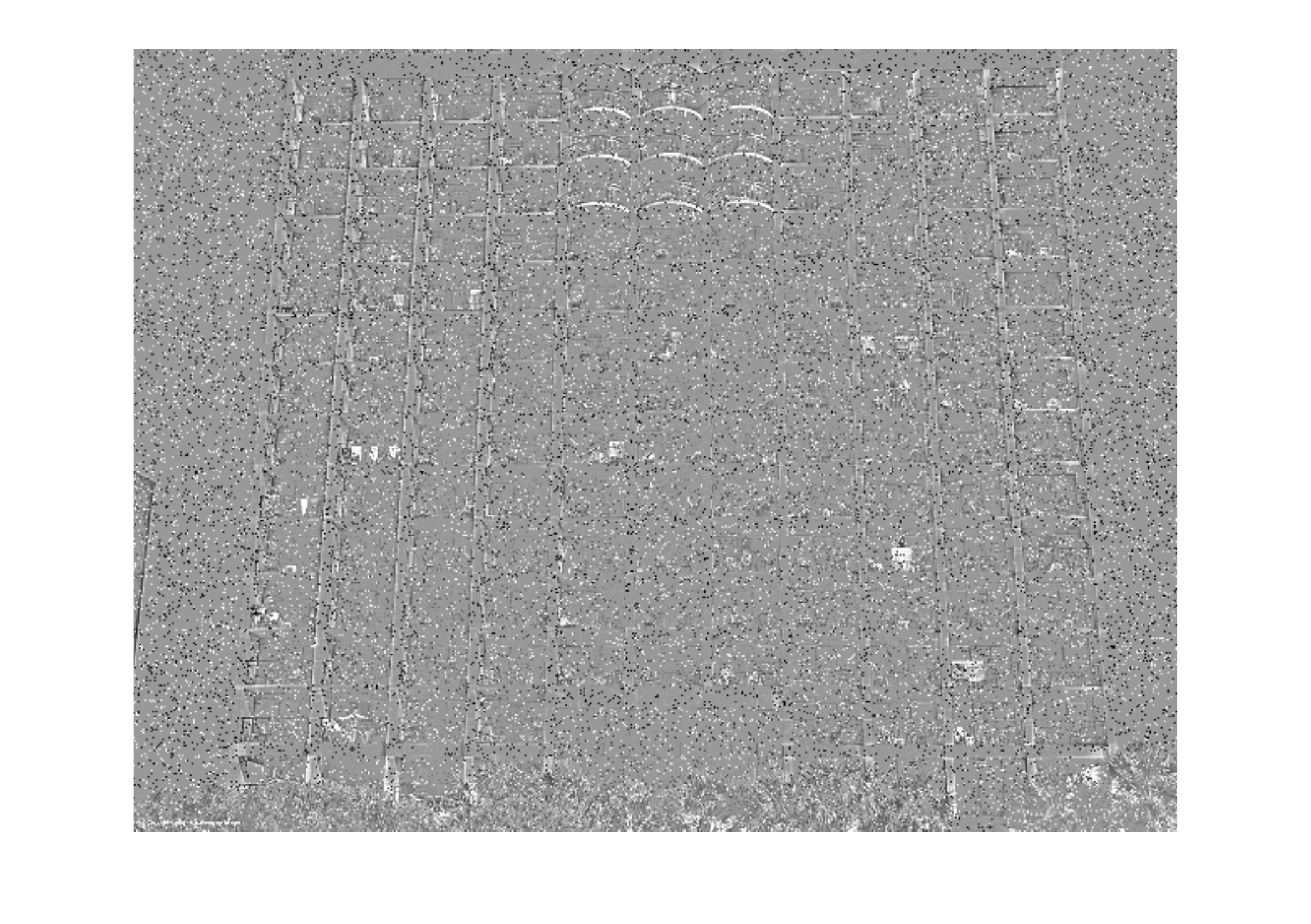}
                \label{fig2:sparse}}
\subfloat[smooth component]{
        \includegraphics[trim=1cm 1.2cm 1cm 8mm, clip,width=0.45\textwidth]{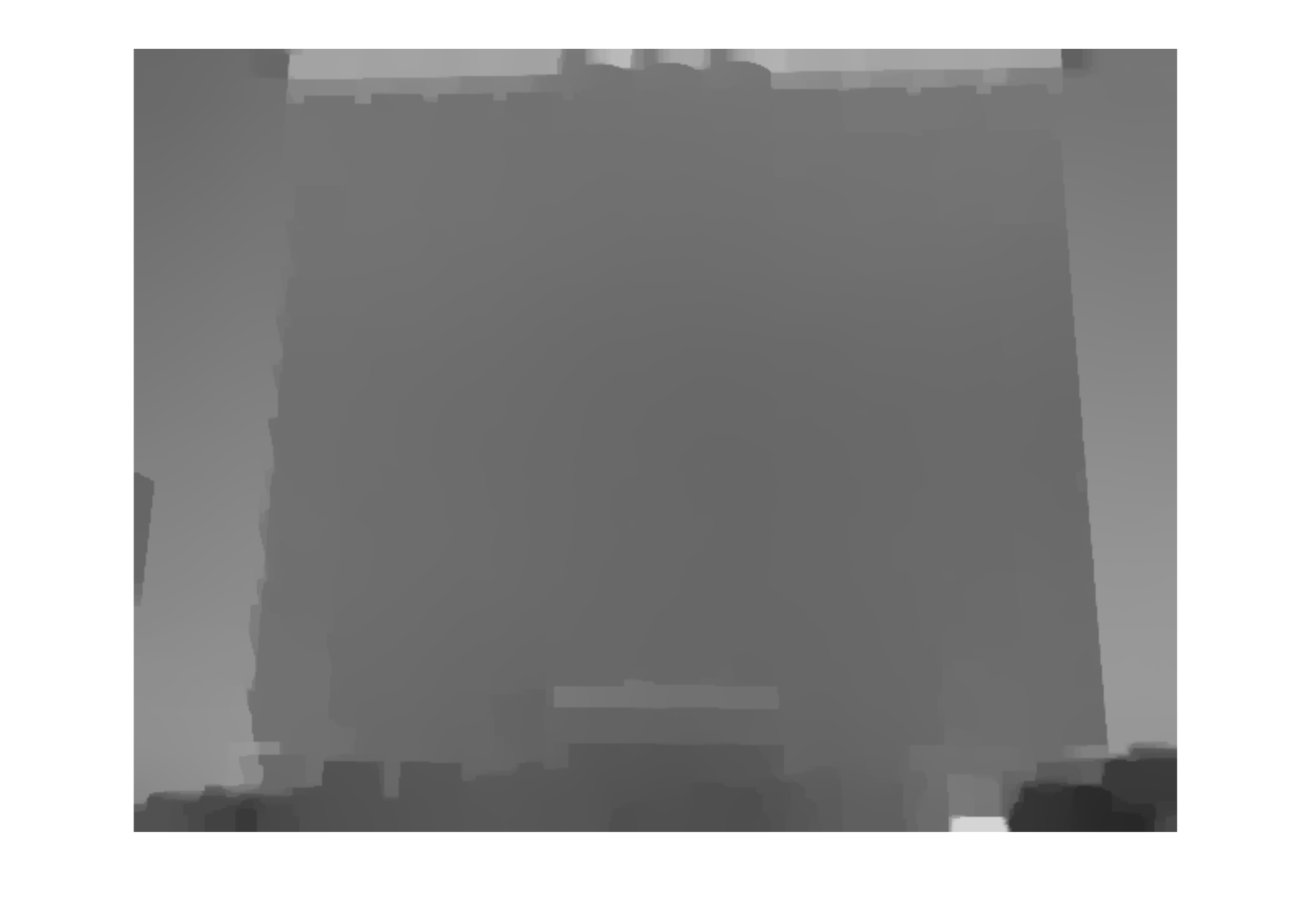}
                \label{fig2:smooth}}
\caption{Observed and decomposed images (size $480\times640$)}
\label{building}
\end{figure}


\section{Semi-Separable Convex Problems}\label{AltDir}
\subsection{Preliminaries}\label{section: Preliminaries}
Our problem of interest in this section is  problem {(\ref{model:constrained:0}), (\ref{model:saddlepoint1}),} namely,
\be
\begin{array}{rl}
\Opt&=\min\limits_{[y^1;\ldots;y^K] \in Y_1\times \cdots\times Y_K}
\left\{f([y^1;\ldots;y^K] ):= \sum_{k=1}^K[\psi_k(y^k)+\Psi_k(y^k)]: \;\sum_{k=1}^K A_ky^k =b\right\}\\
&=\min\limits_{[y^1;\ldots;y^K] \in Y_1\times \cdots\times Y_K}\left\{
\sum\limits_{k=1}^K\left[\psi_k(y^k)+\Psi_k(y^k)\right]:\;g([y^1;...;y^K])\leq 0 \right\},\\
&\qquad g([y^1;...;y^K])=\pi^*\left(\sum\limits_{k=1}^K A_ky^k-b\right)=\max\limits_{\pi(w)\le 1} \sum\limits_{k=1}^K\langle A_ky^k-b,w
\rangle,
\end{array}
\ee{model:constrained:00}
where $\pi(\cdot)$ is some norm and $\pi^*(\cdot)$ is the conjugate norm.  A straightforward approach to \rf{model:constrained:00} would be to rewrite it as a saddle point problem
\be
\min_{[y^1;\ldots;y^K] \in Y_1\times \cdots\times Y_K}\max_{w}\left\{\sum_{k=1}^K[\psi_k(y^k)+\Psi_k(y^k)] +\langle \sum_{k=1}^KA_kz^k-b,w\rangle\right\}
\ee{unboundedw}
and solve by the mirror-prox algorithm from section \ref{MP} adjusted to work with an unbounded domain $U$, or, alternatively, we could replace $\max_w$ with $\max_{w:\;\pi(w)\leq R}$ with ``large enough'' $R$ and use the above algorithm ``as is.'' The potential problem with this approach is that if the {$w$-component $w^*$ of the saddle point of} \rf{unboundedw} is of large $\pi$-norm (or ``large enough'' $R$ is indeed large), the (theoretical) efficiency estimate would be bad since it is proportional to the magnitude of $w^*$ (resp., to $R$). To circumvent this difficulty, we apply to \rf{model:constrained:00} the sophisticated policy originating from \cite{LNN}. {This policy requires the set $Y=Y_1\times...\times Y_K$ to be bounded, which we assume below.}

\paragraph{Course of actions.}
Note that our problem of interest is of the generic form
\begin{equation}\label{octeq34}
\Opt=\min_{y\in Y}\left\{f(y): \;g(y)\leq0\right\}
\end{equation}
where $Y$ is a convex {compact} set in a Euclidean space $E$, {$f$  and $g:\;Y\to\bR$ are  convex} and Lipschitz continuous functions.
For the time being, we focus on (\ref{octeq34}) and assume that the problem is {feasible and thus} solvable.
\par
We intend to solve (\ref{octeq34}) by the generic algorithm presented in \cite{LNN}; for our now purposes, the following description of the algorithm will do:
\begin{enumerate}
\item The algorithm works in {\sl stages}. Stage $s=1,2,...$ is associated with {\sl working parameter} $\alpha_s\in(0,1)$. We set $\alpha_1=\half$.
\item At stage $s$, we apply a first order method $\cB$ to the problem
 \be
\hskip-80pt{(P_s)}~~~~~~~~~~~~~~~~~~~~~~~
\Opt_s=\min_{y\in Y} \left\{f_s(y)=\alpha_s f(y)+(1-\alpha_s) g(y)\right\}
\ee{P_s}
The only property of the algorithm $\cB$ which matters here is its ability, when run on $(P_s)$, to produce in course of $t=1,2,...$ steps iterates $y_{s,t}$, upper bounds $\overline{f}_s^t$ on $\Opt_s$ and lower bounds $\underline{f}_{s,t}$ on $\Opt_s$ in such a way that
\begin{enumerate}
\item for every $t=1,2,...$, the $t$-th iterate $y_{s,t}$ of $\cB$ as applied to $(P_s)$ belongs to $Y$;
\item the upper bounds $\overline{f}_s^t$ are nonincreasing in {$t$} (this is ``for free'')  and ``are achievable,'' that is, they are of the form
$$
\overline{f}_s^t=f_s(y^{s,t}),
$$
where $y^{s,t}\in Y$ is a vector which we have at our disposal at step $t$ of stage $s$;
\item the lower bounds $\underline{f}_{s,t}$ should be nondecreasing in {t} (this again is ``for free'');
\item for some nonincreasing sequence $\epsilon_t\to+0$, $t\to\infty$, we should have
$$
\overline{f}_s^t-\underline{f}_{s,t}\leq\epsilon_t
$$
for all $t$ and $s$.
\end{enumerate}
Note that since (\ref{octeq34}) is solvable, we clearly have $\Opt_s\leq \alpha_s\Opt$, implying that the quantity $\underline{f}_{s,t}/\alpha_s$ is a lower bound on $\Opt$. Thus, at step $t$ of stage $s$ we have at our disposal a number of valid lower bounds on $\Opt$; we denote the best (the largest) of these bounds $\underline{\Opt}_{s,t}$, so that
\begin{equation}\label{octeq36}
\Opt\geq \underline{\Opt}_{s,t}\geq \underline{f}_{s,t}/\alpha_s
\end{equation}
for all $s,t$, and $\underline{\Opt}_{s,t}$  is nondecreasing in time\footnote{in what follows, we call a collection $a_{s,t}$ of reals nonincreasing in time, if $a_{s',t'}\leq a_{s,t}$ whenever $s'\geq s$, same as whenever $s=s'$ and $t'\geq t$. ``Nondecreasing in time'' is defined similarly.}.
\item When the First Order oracle is invoked at step $t$ of stage $s$, we get at our disposal a triple $(y_{s,t}\in Y,f(y_{s,t}),g(y_{s,t}))$. We assume that all these triples are somehow memorized. Thus, after calling First Order oracle at step $t$ of stage $s$, we have at our disposal a finite set $Q_{s,t}$ on the 2D plane such that for every point $(p,q)\in Q_{s,t}$ we have at our disposal a vector $y_{pq}\in Y$ such that
$f(y_{pq})\leq p$ and $g(y_{pq})\leq q$; the set $Q_{s,t}$ (in today terminology, {\sl a filter}) is comprised of all pairs $(f(y_{s',t'}),g(y_{s',t'}))$ generated so far. We set
\begin{equation}\label{octeq37}
\begin{array}{rcl}
h_{s,t}(\alpha)&=&\min_{(p,q)\in Q_{s,t}} \left[\alpha (p-\underline{\Opt}_{s,t}) + (1-\alpha) q\right]: [0,1]\to \bR,\\
\Gap(s,t)&=&\max\limits_{0\leq\alpha\leq 1} h_{s,t}(\alpha).\\
\end{array}
\end{equation}
\item  Let $\Delta_{s,t}=\{\alpha\in[0,1]: h_{s,t}(\alpha)\geq0\}$, so that $\Delta_{s,t}$ is a segment in $[0,1]$. Unless we have arrived at $\Gap(s,t)=0$ (i.e., got an optimal solution to (\ref{octeq34}), see (\ref{octeq38})), $\Delta_{s,t}$ is not a singleton (since otherwise $\Gap(s,t)$ were 0). Observe also that $\Delta_{s,t}$ are nested: $\Delta_{s',t'}\subset \Delta_{s,t}$ whenever $s'\geq s$, same as whenever $s'=s$ and $t'\geq t$.
    \par
We continue iterations of stage $s$ while  $\alpha_s$ is ``well-centered'' in $\Delta_{s,t}$, e.g., belongs to the mid-third of the segment. When this condition is violated, we start stage $s+1$, specifying $\alpha_{s+1}$ as the midpoint of $\Delta_{s,t}$.
\end{enumerate}
The properties of the aforementioned   routine are summarized in the following statement (cf. \cite{LNN}).
\begin{proposition}\label{pro:LNN-filter}
\item{\rm (i)}
$\Gap(s,t)$ is {nonincreasing in time}. Furthermore, at step $t$ of stage $s$, we have at our disposal a solution $\widehat{y}^{s,t}\in Y$ to {\rm\rf{octeq34}} such that
\begin{equation}\label{octeq38}
f(\widehat{y}^{s,t})\leq\Opt + \Gap(s,t),\;\;\mbox{and}\;\; g(\widehat{y}^{s,t}) \leq \Gap(s,t),
\end{equation}
so that $\widehat{y}^{s,t}$ belongs to the domain $Y$ of problem (\ref{octeq34}) and is both $\Gap(s,t)$-feasible and $\Gap(s,t)$-optimal.
\item{\rm (ii)} For every $\epsilon>0$, the number $s(\epsilon)$ of stages until a pair $(s,t)$ with $\Gap(s,t)\leq\epsilon$ is found obeys the bound
\begin{equation}\label{octeq39}
s(\epsilon)\leq{\ln\left({3L\epsilon^{-1}}\right)\over \ln\left(4/3\right)},
\end{equation}
where $L<\infty$ is an a priori upper bound on $\max_{y\in Y}\max[|f(y)|,|g(y)|]$. Besides this,
the number of steps at each stage does not exceed
\begin{equation}\label{octeq40}
T(\epsilon)=\min\{t\geq 1: \epsilon_t\leq {\epsilon\over 3}\}+1.
\end{equation}
\end{proposition}

\subsection{Composite Mirror Prox algorithm for Semi-Separable Optimization}\label{section:our algorithm}
We are {about} to apply the approach above to the semi-separable problem {(\ref{model:constrained:00}), (\ref{model:saddlepoint1})}.
\paragraph{Problem setup} we consider now is as follows (cf. section \ref{sec:41}). For every $k$, $1\leq k\leq K$, we are given
\begin{enumerate}
\item 
Euclidean spaces $E_k$ and $\oE_k$ along with their nonempty closed {and bounded} convex subsets $Y_k$ and {$Z_k$, respectively;}
\item 
proximal setups for $(E_k,Y_k)$ and $(\oE_k, Z_k)$, that is, norms $p_k(\cdot)$ on $E_k$, norms $q_k$ on $\oE_k$, and {d.g.f.'s} $\omega_k(\cdot):Y_k\to \bR$, $\oomega_k(\cdot):Z_k\to \bR$, which are {compatible with} $p_k(\cdot)$ {and} $q_k(\cdot)$, {respectively;}
\item 
linear mapping $y^k\mapsto A_k y^k:E_k\to E$, where $E$ is a Euclidean space;
\item 
Lipschitz continuous convex functions $\psi_k(y^k):Y_k\to\bR$ along with their \emph{saddle point representations}
\be
\psi_k(y^k)=\sup_{z^k\in Z_k}{[}\phi_k(y^k,z^k)-\oPsi_k(z^k){]},\;\;1\le k\le K,
\ee{model:saddlepoint}
where $\phi_k(y^k,z^k):Y_k\times Z_k\to\bR$ are smooth (with Lipschitz continuous gradients)  functions convex in $y^k\in Y_k$ and concave in $z^k\in Z_k$, and  $\oPsi_k(z^k):Z_k\to\bR$ are Lipschitz continuous convex functions  such that the problems of the form
\begin{equation}
\begin{array}{rl}\min\limits_{z^k\in Z_k}&\left[\oomega_k(z^k)+\langle \xi^k,z^k\rangle +\alpha\oPsi_k(z^k)\right]\quad[\alpha>0]
\end{array}
\end{equation}
are easy to solve;
\item 
Lipschitz continuous convex functions $\Psi_k(y^k):Y_k\to\bR$ such that the problems of the form
\[
\begin{array}{rl}\min\limits_{y^k\in Y_k}&\left[\omega_k(y^k)+\langle \xi^k,y^k\rangle +\alpha\Psi_k(y^k)\right]\quad[\alpha>0]
\end{array}
\]
are easy to solve;
\item  {a} norm $\pi^*(\cdot)$ on $E$, with conjugate norm $\pi(\cdot)$, along with a {d.g.f.} $\widehat{\omega}(\cdot):\;W:=\{w\in E:\pi(w)\leq1\}\to\bR$ {compatible with} $\pi(\cdot)$ and is such that  problems {of the form}
\[
\min_{w\in W}\left[\widehat{\omega}(w)+\langle \xi,w\rangle\right]
\]
are easy to solve.
\end{enumerate}
The outlined data define the sets
\[
\begin{array}{rcl}
Y^+_k&=&\{[y^k;\tau^k]: \;y^k\in Y_k,\tau^k\geq\Psi_k(y^k)\}\subset E_k^+:=E_k\times \bR,\,\,1\leq k\leq K,\\
Z^+_k&=&\{[z^k;\sigma^k]: \;z^k\in Z_k,\sigma^k\geq\oPsi_k(z^k)\}\subset \oE_k^+:=\oE_k\times \bR,\,\,1\leq k\leq K.
\end{array}
\]
The problem of interest here is problem {(\ref{model:constrained:00}), (\ref{model:saddlepoint}):}
\be
&\begin{array}{rl}
\Opt=\min\limits_{[y^1;\ldots;y^K]}\max\limits_{[z^1;\ldots;z^K]}
&\Big\{ \sum_{k=1}^K[\phi_k(y^k,z^k)+\Psi_k(y^k)-\oPsi_k(z^k)]: \;\;\pi^*\left(\sum\limits_{k=1}^K A_ky^k-b\right)\le 0,\\
&[y^1;\ldots;y^K] \in Y_1\times \cdots\times Y_K,\;\;[z^1;\ldots;z^k]\in Z_1\times \cdots\times Z_K
\Big\} \end{array}\nn
&\begin{array}{rl}
=\min\limits_{\{[y^k;\tau^k]\}_{k=1}^K}\max\limits_{\{[z^k;\sigma^k]\}_{k=1}^K}
&\Big\{ \sum_{k=1}^K[\phi_k(y^k,z^k)+\tau^k-\sigma^k]:\;\;
\max\limits_{w\in W} \sum\limits_{k=1}^K\langle A_ky^k-b,w
\rangle\leq 0,\\
&\{[y^k;\tau^k]\in Y^+_k\}_{k=1}^K,\;\;\{[z^k;\sigma^k]\in Z^+_k\}_{k=1}^K,\;w\in W
\Big\}. \end{array}
\ee{model:constrained}
\paragraph{Solving  \rf{model:constrained}} using the approach in the previous section amounts to resolving a sequence of problems $(P_s)$ as in \rf{P_s} where, with a slight abuse of notation,
\bse
Y&=&\left\{y=\{[y^k;\tau^k]\}_{k=1}^K:\;[y^k;\tau^k]\in Y^+_k, \;\tau^k\leq C_k,\,1\leq k\leq K\right\};\\
f(y)&=&\max\limits_{z=\{[z^k;\sigma^k]\}_{k=1}^K}\left\{\sum_{k=1}^K[\phi_k(y^k,z^k)+\tau^k-\sigma^k]: \;z\in Z=\{[z^k;\sigma^k]\in Z^+_k\}_{k=1}^K\right\};\\
g(y)&=&\max\limits_{w}\left\{\sum_{k=1}^K\langle A_ky^k-b,w\rangle:\;w\in W\right\}.\\
\ese
Here $C_k\geq\max_{y^k\in Y_k}\Psi_k(y^k)$ are finite constants introduced to make $Y$ compact, as required in the premise of Proposition \ref{pro:LNN-filter}; it is immediately seen that the magnitudes of these constants (same as their very presence) does not affect the algorithm $\cB$ we are about to {describe}.\par
The algorithm $\cB$ we intend to use will solve $(P_s)$ by reducing the problem to the saddle point problem
\be
\begin{array}{rl}\overline{\Opt}=\min\limits_{y}\;\max\limits_{[z;w]} &\Big\{
\Phi(y,[z;w]):=\alpha\sum_{k=1}^K[\phi_k(y^k,z^k)+\tau^k-\sigma^k]+(1-\alpha)\sum_{k=1}^K\langle A_ky^k-b,w\rangle:\\
&
\multicolumn{1}{r}{y=\{[y^k;\tau^k]\}_{k=1}^K\in Y,\;[z=\{[z^k;\sigma^k]\}_{k=1}^K\in Z;\; w\in W]\Big\},}\\
\end{array}\ee{octeq42}
where $\alpha=\alpha_s$.
\par
Setting
\bse
U&=&\{u=[y^1;...;y^K;z^1;...;z^K;w]:\;y^k\in Y_k,\;z^k\in Z_k,\;1\leq k\leq K, w\in W\},\\
X&=&\{[u;v=[\tau^1;...;\tau^K;\sigma^1;...\sigma^K]]:\; u\in U, \;\Psi_k(y^k)\leq \tau^k\leq C_k,\;\oPsi_k(z^k)\le \sigma^k,\;1\leq k\leq K\},
\ese
$X$ can be thought of as the domain of the variational inequality associated with (\ref{octeq42}), the monotone operator in question being
\begin{equation}\label{octeq43}
\begin{array}{rcl}
F(u,v)&=&[F_u(u);F_v],\\
F_u(u)&=&\left[\begin{array}{l}
\left\{
\alpha\nabla_y\phi_k(y^k,z^k)+(1-\alpha)A_k^Tw\right\}_{k=1}^K\\
\left\{
-\alpha\nabla_z\phi_k(y^k,z^k)\right\}_{k=1}^K\\
(1-\alpha)[b-\sum_{k=1}^KA_ky^k]
\end{array}\right],\\
F_v&=&\alpha[1;...;1].\\
\end{array}
\end{equation}
By exactly the same reasons as in  section \ref{MultiTerm}, with properly assembled norm on the embedding space of $U$ and d.g.f., (\ref{octeq42}) can be solved by the Mirror Prox algorithm from section \ref{MP}. Let us denote
$$
\zeta^{s,t}=
\left[\widehat{y}^{s,t}=\{[\widehat{y}^k;\widehat{\tau}^k]\}_{k=1}^K\in Y;\left[z^{s,t}\in Z; w^{s,t}\in W\right]\right]
$$
the approximate solution obtained in course of $t=1,2,...$ steps of CoMP when solving $(P_s)$, and let
\[
\widehat{f}_s^t:=\max_{z\in Z,w\in W}\Phi(\widehat{y}^{s,t},[z;w])=\alpha\sum_{k=1}^K[\psi_k(\widehat{y}^k)+\widehat{\tau}^k]+(1-\alpha)\pi^*\left(\sum_{k=1}^KA_k\widehat{y}^k-b\right)
\]
be the corresponding value of the objective of $(P_s)$. It holds
\begin{equation}\label{octeq44}
\widehat{f}_s^t-\overline{\Opt}\leq \epsilonsad(\zeta^{s,t}\big| Y,Z\times W,\Phi)\leq \epsilon_t:=O(1)\cL/t,
\end{equation}
where $\cL<\infty$ is explicitly given by the proximal setup we use and by the related Lipschitz constant of $F_u(\cdot)$ (note that this constant can be chosen to be independent of $\alpha\in[0,1]$).
We assume that computing the corresponding objective value is a part of step $t$ (these computations increase the complexity of a step by factor at most $O(1)$), and thus that $\overline{f}_s^t\leq \widehat{f}_s^t$. By (\ref{octeq44}), the quantity $\widehat{f}_s^t-\epsilon_t$ is a valid lower bound on the optimal value of $(P_s)$, and thus we can ensure that $\underline{f}_{s,t}\geq \widehat{f}_s^t-\epsilon_t$. The bottom line is that with the outlined implementation, we have
$$
\overline{f}_s^t-\underline{f}_{s,t}\leq\epsilon_t
$$
for all $s,t$, with $\epsilon_t$ given by (\ref{octeq44}). Consequently, by Proposition \ref{pro:LNN-filter}, the total number of CoMP steps needed to find a belonging to the domain of the problem of interest \rf{model:constrained:00}  $\epsilon$-feasible and $\epsilon$-optimal solution to this problem can be upper-bounded by
$$
O(1)\ln\left({3L\over\epsilon}\right) \left( {\cL\over\epsilon}\right),
$$
where $L$ and $\cL$ are readily given by the smoothness parameters of $\phi_k$ and {by} the proximal setup we use.

\subsection{Numerical illustration: $\ell_1$-minimization}
\paragraph{Problem of interest.} We consider the simple $\ell_1$ minimization problem
\begin{equation}\label{model: l1-minimization}
\min\limits_{x\in X} \left\{\|x\|_1:\; Ax = b\right\}
\end{equation}
where $x\in\bR^n${,} $A\in\bR^{m\times n}$ and $m<n$. {Note that} this problem can also be written in the semi-separable form
\begin{equation*}
\begin{array}{c}
\min\limits_{x\in X} \left\{\sum_{k=1}^K\|x_k\|_1: \;\sum_{k=1}^KA_kx_k = b\right\}
\end{array}
\end{equation*}
if the data is partitioned into $K$ blocks: $x=[x_1;x_2;\ldots;x_K]$ and $A=[A_1,A_2,\ldots,A_K]$.

Our main purpose here is to test the approach described in \ref{section: Preliminaries}  and compare it to the {simplest} approach {where we directly apply CoMP to the (saddle point reformulation of the) problem $\min_{x\in X}[\|x\|_1+R\|Ax-b\|_2]$ with large enough value of $R$.} {For} the sake of simplicity, {we work with} the case when $K=1$ and $X=\{x\in\bR^n:\|x\|_2\leq 1\}$.

\paragraph{Generating the data.}  In the experiments to be reported, the data of (\ref{model: l1-minimization}) were generated as follows. Given $m,n$, we first build a sparse solution $x^*$ by drawing  random vector from {the} standard Gaussian distribution $\cN(0,I_n)$, zeroing out part of the entries and {scaling the resulting vector to enforce $x^*\in X$. We also} build a dual solution $\lambda^*$ by scaling a random vector from distribution $\cN(0,I_m)$ to satisfy $\|\lambda^*\|_2=R_*$ for a prescribed $R_*$. Next we generate $A$ and $b$ such that $x^*$ and $\lambda^*$ are indeed the optimal primal and dual solutions to the $\ell_1$ minimization problem} (\ref{model: l1-minimization}), i.e. $A^T\lambda^*\in\partial\big|_{x=x^*}\|x\|_1$ and $Ax^*=b$. To achieve this, we set
\begin{equation*}
\begin{array}{c}
A = \frac{1}{\sqrt{n}}\hat{F}_n+ pq^T,\; b = Ax^*
\end{array}
\end{equation*}
where $ p=\frac{\lambda^*}{\|\lambda^*\|_2^2},  q\in\partial\big|_{x=x^*} \|x\|_1-\frac{1}{\sqrt{n}}\hat{F}_n \lambda^*$, and  $\hat{F}_n$ is a $m\times n$ submatrix randomly selected from the DFT matrix $F_n$.
We expect that the larger is the $\|\cdot\|_2$-norm $R_*$ of the dual solution, the harder is  problem (\ref{model: l1-minimization}).
\paragraph{Implementing the algorithm.} When implementing the algorithm {from} section \ref{section:our algorithm}, we apply {at each stage $s=1,2,...$} CoMP to the saddle point problem
\[
(P_s): \quad \min_{\substack{x,\tau:\; \|x\|_2\leq 1,  \tau\geq\|x\|_1}}\max_{w:\|w\|_2\leq 1} \left\{\alpha_s \tau +(1-\alpha_s)\langle Ax-b,w\rangle\right\}.
\]
{The proximal setup for CoMP is given by equipping the embedding space of $U=\{u=[x;w]:x\in X,\|w\|_2\leq1\}$ with the norm $\|u\|_2=\sqrt{\frac{1}{2}\|x\|_2^2+\frac{1}{2}\|w\|_2^2}$ and equipping $U$ with the d.g.f. $\omega(u)=\frac{1}{2}\|x\|_2^2+\frac{1}{2}\|w\|_2^2$.} In the {sequel} we refer {to the resulting algorithm} as {\sl sequential} CoMP.  {For comparison, we solve the same problem by applying CoMP to the saddle point problem}
\[
(P_R): \quad  \min_{\substack{x,\tau:\;\|x\|_2\leq 1, \tau\geq\|x\|_1}}\max_{w:\|w\|_2\leq 1} \left\{\tau+R\langle Ax-b,w\rangle\right\}
\]
with $R=R_*$; the resulting algorithm is referred to as {\sl simple} CoMP. Both sequential CoMP and simple  CoMP algorithms  are terminated when the relative nonoptimality and constraint violation are both less than  $\epsilon=10^{-5}$, namely,
\begin{equation*}
\begin{array}{c}
\epsilon(x) :=\max\left\{\frac{\|x\|_1-\|x_*\|_1}{\|x_*\|_1}, \|Ax-b\|_2\right\}\leq 10^{-5}.
\end{array}
\end{equation*}

\begin{table}[h!]
\begin{center}{\scriptsize
 \begin{tabular}{|ccc||cr||cr||}
 \hline
$n$&$m$&$c$&\multicolumn{2}{c||}{\textbf{sequential CoMP}}&\multicolumn{2}{c||}{\textbf{simple CoMP}}\\
 \cline{4-7}
&&$(R_*=c\cdot n)$ & steps & CPU(sec) & steps & CPU(sec) \\
\hline
1024	&512	&1	&7653	&18.68	&31645	&67.78\\
&		                   &5	&43130	&44.66	&90736	&90.67\\
&			 		&10	&48290	&49.04	&93989	&93.28\\
\hline							
4096	&2048	&1	&28408	&85.83	&46258	&141.10\\
&					&5	&45825	&199.96	&93483	&387.88\\
&					&10	&52082	&179.10	&98222	&328.31\\
\hline							
16384	&8192	&1  &43646	&358.26	&92441	&815.97\\
&					&5	&48660	&454.70	&93035	&784.05\\
&					&10&55898	&646.36	&101881&1405.80\\
\hline							
65536	&32768	&1	&45153	&3976.51	&92036	&4522.43\\
&			             &5  &55684	&4138.62  &100341 &8054.35\\
&					&10&69745	&6214.18	&109551&9441.46\\
\hline
262144&131072	&1	&46418	&6872.64	&96044	&14456.99\\
&						&5	&69638	&10186.51	&109735&16483.62\\
&						&10	&82365	&12395.67	&95756	&13634.60\\
\hline
 \end{tabular}}
\end{center}
\caption{\label{table: l1-minimization}$\ell_1$-minimization. Platform: ISyE Condor Cluster}
 \end{table}

\paragraph{Numerical results} are presented in Table \ref{table: l1-minimization}. One can immediately see that {to achieve the desired accuracy, the simple CoMP} with $R$ set to  $R_*$, i.e., to the exact magnitude of the true Lagrangian multiplier, requires almost twice as many steps as the sequential CoMP. In more realistic examples, the simple CoMP will additionally suffer from the fact that the magnitude of the optimal Lagrange multiplier is not known in advance, and the penalty $R$ in $(P_R)$ should be somehow tuned ``online.''

%
%
\bibliographystyle{abbrv}  
\bibliography{MyLibrary}

\appendix
\section{Proof of Theorem \ref{theMP}}\label{sect:prooftheMP}
\paragraph{0$^o$.} Let us verify that the prox-mapping (\ref{eq:prox}) indeed is well defined whenever $\zeta=\gamma F_v$ with $\gamma>0$. All we need is to show that whenever $u\in U$, $\eta\in E_u$, $\gamma>0$ and $[w_t;s_t]\in X$, $t=1,2,...$, are such that $\|w_t\|_2+\|s_t\|_2\to\infty$ as $t\to\infty$, we have
$$
r_t:=\underbrace{\langle \eta-\omega'(u),w_t\rangle+\omega(w_t)}_{a_t}+\underbrace{\gamma\langle F_v,s_t\rangle}_{b_t} \to\infty ,\,t\to\infty.
$$
Indeed, assuming the opposite and passing to a subsequence, we make the sequence $r_t$ bounded. Since $\omega(\cdot)$ is strongly convex, modulus 1, w.r.t. $\|\cdot\|$, and the linear function $\langle F_v,s\rangle$ of $[w;s]$ is below bounded on $X$ by {\bf A4},  boundedness of the sequence $\{r_t\}$ implies
boundedness of the sequence $\{w_t\}$, and since $\|[w_t;s_t]\|_2\to\infty$ as $t\to\infty$, we get $\|s_t\|_2\to\infty$ as $t\to\infty$. Since $\langle F_v,s\rangle$ is coercive in $s$ on $X$ by {\bf A4}, and $\gamma>0$, we conclude that $b_t\to\infty$, $t\to\infty$, while the sequence $\{a_t\}$ is bounded since the sequence $\{w_t\in U\}$ is so and $\omega$ is continuously differentiable. Thus, $\{a_t\}$ is bounded, $b_t\to\infty$, $t\to\infty$, implying that $r_t\to\infty$, $t\to\infty$, which is the desired contradiction

\paragraph{1$^o$}. Recall  the well-known identity \cite{CT93}: for all $u,u',w\in U$ one has
\begin{equation}\label{threetermid}
\langle V'_{u}(u'),w-u'\rangle =V_{u}(w)-V_{u'}(w)-V_{u}(u').
\end{equation}
{
Indeed, the right hand side is
\bse
\lefteqn{[\omega(w)-\omega(u)-\langle\omega'(u),w-u\rangle]-[\omega(w)-\omega(u')-\langle\omega'(u'),w-u'\rangle]-[\omega(u')-\omega(u)-\langle\omega'(u),u'-u\rangle]}\\
&=&\langle \omega'(u),u-w\rangle +\langle \omega'(u),u'-u\rangle +\langle \omega'(u'),w-u'\rangle =\langle\omega'(u')- \omega'(u),w-u'\rangle=\langle V'_{u}(u'),w-u'\rangle.
\ese}
For $x=[u;v]\in X,\;\xi=[\eta;\zeta]$, let $P_x(\xi)=[u';v']\in X$. By the optimality condition for the problem \rf{eq:prox}, for all $[s;w]\in X$,
\[
\langle \eta+V'_u(u'),u'-s\rangle+\langle\zeta,v'-w\rangle \le 0,
\]
which by \rf{threetermid} implies that
\be
\langle \eta,u'-s\rangle+\langle\zeta,v'-w\rangle \le \langle V'_u(u'),s-u'\rangle=
V_{u}(s)-V_{u'}(s)-V_{u}(u').
\ee{prox_lemma}
\paragraph{2$^o$.} When applying \rf{prox_lemma} with $[u;v]=[u_\tau;v_\tau]=x_\tau$, $\xi=\gamma_\tau F(x_\tau)=[\gamma_\tau F_u(u_\tau);\gamma_\tau F_v]$, $[u';v']=[u'_\tau;v'_\tau]=y_\tau$, and $[s;w]=[u_{\tau+1};v_{\tau+1}]=x_{\tau+1}$ we obtain:
\be
\gamma_\tau [\langle F_u(u_\tau),u'_\tau-u_{\tau+1}\rangle+\langle F_v,v'_\tau-v_{\tau+1}\rangle]\le V_{u_\tau}(u_{\tau+1})-V_{u'_\tau}(u_{\tau+1})-V_{u_\tau}(u'_{\tau});
\ee{prox100}
and applying \rf{prox_lemma} with $[u;v]=x_\tau$, $\xi=\gamma_\tau F(y_\tau)$, $[u';v']=x_{\tau+1}$, and $[s;w]=z\in X$ we get:
\be
\gamma_\tau [\langle F_u(u'_\tau),u_{\tau+1}-s\rangle+\langle F_v,v_{\tau+1}-w\rangle]\le V_{u_\tau}(s)-V_{u_{\tau+1}}(s)-V_{u_\tau}(u_{\tau+1}).
\ee{prox101}
Adding \rf{prox101} to \rf{prox100} we obtain for every $z=[s;w]\in X$
\be
\lefteqn{\gamma_\tau\langle F(y_\tau),y_\tau-z\rangle=\gamma_\tau[\langle F_u(u'_\tau),u'_\tau-s\rangle+\langle F_v,v'_\tau-w\rangle]}\nn
&\le & V_{u_\tau}(s)-V_{u_{\tau+1}}(s)+\underbrace{\gamma_\tau\langle F_u(u'_\tau)-F_u(u_\tau),u'_\tau-u_{\tau+1}\rangle-V_{u'_\tau}(u_{\tau+1})-V_{u_\tau}(u'_{\tau})
}_{\delta_\tau}.
\ee{prox102}
Due to the strong convexity, modulus 1, of $V_u(\cdot)$ w.r.t. $\|\cdot\|$, $V_u(u')\geq {1\over 2}\|u-u'\|^2$ for all $u,u'$. Therefore,
\bse
\delta_\tau&\leq& \gamma_\tau\|F_u(u'_\tau)-F_u(u_\tau)\|_*\|u'_\tau-u_{\tau+1}\|-\half\|u'_\tau-u_{\tau+1}\|^2-\half\|u_\tau-u'_\tau\|^2\\
&\leq&
\half\left[\gamma_\tau^2\|F_u(u'_\tau)-F_u(u_\tau)\|_*^2-\|u_\tau-u'_\tau\|^2\right]\\
&\leq&
\half\left[\gamma_\tau^2[M+L\|u'_\tau-u_\tau\|]^2-\|u_\tau-u'_\tau\|^2\right],
\ese
where the last inequality is due to \rf{ML}. Note that $\gamma_\tau L<1$ implies that
$$
\gamma_\tau^2[M+L\|u'_\tau-u_\tau\|]^2-\|u'_\tau-u_\tau\|^2\leq \max_r
 \left[\gamma_\tau^2[M+Lr]^2-r^2\right]={\gamma_\tau^2M^2\over1-\gamma_\tau^2L^2}.
$$
Let us assume that the stepsizes $\gamma_\tau>0$  ensure that \rf{gammaupperbound} holds, meaning that
$
\delta_\tau \leq \gamma_\tau^2M^2$
(which, by the above analysis, is definitely the case when $0<\gamma_\tau\leq {1\over\sqrt{2}L}$; when $M=0$, we can take also $\gamma_\tau\leq {1\over L}$).
When summing up inequalities \rf{prox102} over $\tau=1,2,...,t$ and taking into account that $V_{u_{t+1}}(s)\geq0$, we conclude that for all $z=[s;w]\in X$,
\[
\sum_{\tau=1}^t\lambda^t_\tau\langle F(y_\tau),y_\tau-z\rangle\leq
{V_{u_1}(s) +\sum_{\tau=1}^t\delta_\tau\over\sum_{\tau=1}^t\gamma_\tau}\leq{V_{u_1}(s) +M^2\sum_{\tau=1}^t\gamma_\tau^2\over
\sum_{\tau=1}^t\gamma_\tau},\;\;\;\lambda^t_\tau=\gamma_\tau/\sum_{i=1}^t\gamma_i.\;\;\;\hbox{\qed}
\]

\section{Proof of Lemma \ref{lemnew1}}
{\bf Proof.} All we need to verify is the second inequality in (\ref{thenwehave}). To this end note that when $t=1$, the inequality in (\ref{thenwehave}) holds true by definition of $\widehat{\Theta}(\cdot)$. Now let $1<t\leq N+1$. Summing up the inequalities (\ref{prox102}) over $\tau=1,...,t-1$, we get for every $x=[u;v]\in X$:
$$
\sum_{\tau=1}^{t-1}\langle F(y_\tau),y_\tau-[u;v]\rangle \leq V_{u_1}(u)-V_{u_t}(u)+\sum_{\tau=1}^{t-1}\delta_\tau\leq V_{u_1}(u)-V_{u_t}(u)+\sum_{\tau=1}^{t-1}\delta_\tau\leq
V_{u_1}(u)-V_{u_t}(u)+M^2\sum_{\tau=1}^{t-1}\gamma_\tau^2
$$
(we have used (\ref{gammaupperbound})). When  $[u;v]$ is $z_*$, the left hand side in the resulting inequality is $\geq0$, and we arrive at
$$
V_{u_t}(u_*)\leq V_{u_1}(u_*)+M^2\sum_{\tau=1}^{t-1}\gamma_\tau^2,
$$
whence
$$
{1\over 2}\|u_t-u_*\|^2\leq V_{u_1}(u_*)+M^2\sum_{\tau=1}^{t-1}\gamma_\tau^2
$$
whence also
$$
\|u_t-u_1\|^2\leq 2\|u_t-u_*\|^2+2\|u_*-u_1 \|^2\leq 4[V_{u_1}(u_*)+M^2\sum_{\tau=1}^{t-1}\gamma_\tau^2]+4V_{u_1}(u_*)
$$
and therefore
\begin{equation}\label{Roft}
\|u_t-u_1\|\leq 2\sqrt{2V_{u_1}(u_*)+M^2\sum_{\tau=1}^{t-1}\gamma_t^2}= R_N,
\end{equation}
and (\ref{thenwehave}) follows. \qed

\section{Proof of Proposition \ref{propDifferentWeights}}
{\bf Proof.} From (\ref{prox102}) and (\ref{gammaupperbound}) it follows that
$$
\forall(x=[u;v]\in X,\tau\leq N): \lambda_\tau\langle F(y_\tau),y_\tau-x\rangle \leq {\lambda_\tau\over\gamma_\tau}[V_{u_\tau}(u)-V_{u_{\tau+1}}(u)]+M^2\lambda_\tau\gamma_\tau.
$$
Summing up these inequalities over $\tau=1,...,N$, we get $\forall (x=[u;v]\in X)$:
$$
\begin{array}{l}
\sum\limits_{\tau=1}^N\lambda_\tau\langle F(y_\tau),y_\tau-x\rangle\\
 \leq {\lambda_1\over\gamma_1}[V_{u_1}(u)-V_{u_2}(u)]+{\lambda_2\over\gamma_2}[V_{u_2}(u)-V_{u_3}(u)]+...+
{\lambda_N\over\gamma_N}[V_{u_N}(u)-V_{u_{N+1}}(u)]+M^2\sum\limits_{\tau=1}^N\lambda_\tau\gamma_\tau\\
=\underbrace{{\lambda_1\over\gamma_1}}_{\geq0}V_{u_1}(u)+\underbrace{\left[{\lambda_2\over\gamma_2}-{\lambda_1\over\gamma_1}\right]}_{\geq0}V_{u_2}(u)+...+
\underbrace{\left[{\lambda_N\over\gamma_N}-{\lambda_{N-1}\over\gamma_{N-1}}\right]}_{\geq0}
V_{u_N}(u)-{\lambda_N\over \gamma_N}\underbrace{V_{u_{N+1}}(u)}_{\geq0}
+M^2\sum\limits_{\tau=1}^N\lambda_\tau\gamma_\tau\\
\leq {\lambda_1\over\gamma_1}\widehat{\Theta}(\max[R_N,\|u-u_1\|])+\left[{\lambda_2\over\gamma_2}-{\lambda_1\over\gamma_1}\right]\widehat{\Theta}(\max[R_N,\|u-u_1\|])+...\\
\multicolumn{1}{r}{+\left[{\lambda_N\over\gamma_N}-{\lambda_{N-1}\over\gamma_{N-1}}\right]\widehat{\Theta}(\max[R_N,\|u-u_1\|])+M^2\sum\limits_{\tau=1}^N\lambda_\tau\gamma_\tau,}\\
={\lambda_N\over\gamma_N}\widehat{\Theta}(\max[R_N,\|u-u_1\|])+M^2\sum\limits_{\tau=1}^N\lambda_\tau\gamma_\tau,\\
\end{array}
$$
where the concluding inequality is due to (\ref{thenwehave}), and (\ref{newbound}) follows. \qed

\section{Proof of Proposition \ref{pro:LNN-filter}}
 \paragraph{1$^o$.} $h_{s,t}(\alpha)$ are concave piecewise linear functions on $[0,1]$ which clearly are pointwise nonincreasing in time. As a result, $\Gap(s,t)$ is nonincreasing in time.  Further,  we have
\bse
\Gap(s,t)&=&\max\limits_{\alpha\in[0,1]}\left\{\min\limits_{\lambda}
\sum\limits_{(p,q)\in Q_{s,t}}\lambda_{pq}
[\alpha(p-\underline{\Opt}_{s,t})+(1-\alpha)q]:\; \lambda_{pq}\geq0,\;\sum_{(p,q)\in Q_{s,t}}\lambda_{pq}=1\right\}\\
&=&\max_{\alpha\in[0,1]}\sum\limits_{(p,q)\in Q_{s,t}}\lambda^*_{pq}[\alpha (p-\underline{\Opt}_{s,t})+(1-\alpha)q]\\
&=&\max\left[\sum_{(p,q)\in Q_{s,t}}\lambda^*_{pq}(p -\underline{\Opt}_{s,t}),\sum_{(p,q)\in Q_{s,t}}\lambda^*_{pq}q\right],
\ese
where $\lambda^*_{pq}\geq0$ and sum up to 1. Recalling that for every $(p,q)\in Q_{s,t}$ we have at our disposal $y_{pq}\in Y$ such that $p\geq f(y_{pq})$ and  $q\geq g(y_{pq})$, setting $\widehat{y}^{s,t}=\sum\limits_{(p,q)\in Q_{s,t}}\lambda^*_{pq}y_{pq}$ and invoking convexity of $f,g$, we get
$$
f(\widehat{y}^{s,t})\leq \sum\limits_{(p,q)\in Q_{s,t}}\lambda^*_{pq}p\leq \underline{\Opt}_{s,t}+\Gap(s,t), \,\,g(\widehat{y}^{s,t})\leq \sum\limits_{(p,q)\in Q_{s,t}}\lambda^*_{pq}q\leq\Gap(s,t);
$$
and (\ref{octeq38}) follows, due to $\underline{\Opt}_{s,t}\leq \Opt$.
\paragraph{2$^o$.} We have $\overline{f}_s^t=\alpha_sf(y^{s,t}) +(1-\alpha_s)g(y^{s,t})$ for some $y^{s,t}\in Y$ which we have at our disposal at step $t$, implying that $(\bar{p}=f(y^{s,t}),\bar{q}=g(y^{s,t}))\in Q_{s,t}$. Hence by definition of $h_{s,t}(\cdot)$ it holds
$$
h_{s,t}(\alpha_s)\leq \alpha_s (\bar{p}-\underline{\Opt}_{s,t})+(1-\alpha_s)\bar{q}=\overline{f}_s^t-\alpha_s\underline{\Opt}_{s,t}\leq \overline{f}_s^t-\underline{f}_{s,t},
$$
where the concluding inequality is given by (\ref{octeq36}). Thus, $h_{s,t}(\alpha_s)\leq \overline{f}_s^t-\underline{f}_{s,t} \leq\epsilon_t$. On the other hand,
if stage $s$ does not terminate in course of the first $t$ steps, $\alpha_s$ is well-centered in the segment $\Delta_{s,t}$ where the concave function $h_{s,t}(\alpha)$ is nonnegative. We conclude that  $0\leq \Gap(s,t)=\max_{0\leq\alpha\leq 1}h_{s,t}(\alpha)=\max_{\alpha\in\Delta_{s,t}}h_{s,t}(\alpha)\leq 3h_{s,t}(\alpha_s)$. Thus, if a stage $s$ does not terminate in course of the first $t$ steps, we have $\Gap(s,t)\leq 3\epsilon_t$, which implies (\ref{octeq40}). Further, $\alpha_s$ is the midpoint of the segment $\Delta^{s-1}=\Delta_{s-1,t_{s-1}}$, where $t_r$ is the last step of stage $r$  (when $s=1$, we should define $\Delta^0$ as $[0,1]$), and $\alpha_s$ is not well-centered in the segment $\Delta^s=\Delta_{s,t_s}\subset\Delta_{s-1,t_{s-1}}$, which clearly implies that $|\Delta^s|\leq{\rm \small{3\over 4}}|\Delta^{s-1}|$.
Thus, $|\Delta^s|\leq \left({\rm \small{3\over 4}}\right)^s$ for all $s$. On the other hand, when $|\Delta_{s,t}|<1$, we have $\Gap(s,t)=\max_{\alpha\in\Delta_{s,t}}h_{s,t}(\alpha)\leq 3L |\Delta_{s,t}|$ (since $h_{s,t}(\cdot)$ is Lipschitz continuous with constant $3L$ \footnote{we assume w.l.o.g. that $|\underline{\Opt}_{s,t}|\leq L$} and $h_{s,t}(\cdot)$ vanishes at (at least) one endpoint of $\Delta_{s,t}$). Thus, the number of stages before $\Gap(s,t)\leq\epsilon$ is reached indeed obeys the bound (\ref{octeq39}). \qed
\end{document}

%% file: MyMacros.tex
\def\cI{{\cal I}}
\def\SadVal{\hbox{\rm SadVal}}
\def\epsilonsad{{\epsilon_{\hbox{\tiny\rm Sad}}}}
\def\mybox\blacksquare
\def\Bbb#1{{\bf #1}}
\def\Gap{{\hbox{\rm Gap}}}

\def\fnote#1{\footnote}
\def\blacksquare{\hbox{\vrule width 4pt height 4pt depth 0pt}}
\def\square{\hbox{\vrule\vbox{\hrule\phantom{o}\hrule}\vrule}}

\def\nr{\par \noindent}

\def\beq{\begin{equation}}
\def\eeq{\end{equation}}

\newtheorem{theorem}{Theorem}[section]
\newtheorem{lemma}{Lemma}[section]
\newtheorem{corollary}{Corollary}[section]

\newtheorem{proposition}{Proposition}[section]

\newtheorem{remark}{Remark}[section]

\newcommand{\qed}{\hfill $\square$ \nr \medskip}

\def\Tr{\mathop{{\rm Tr}}}

\def\Argmin{\mathop{\rm Argmin}}

\def\Tr{\mathop{{\rm Tr}}}

\def\Argmin{\mathop{\rm Argmin}}

\def\hat{\widehat}

\def\Bbb#1{{\bf #1}}

\def\fnote#1{\footnote}
\def\blacksquare{\hbox{\vrule width 4pt height 4pt depth 0pt}}

\def\square{\hbox{\vrule\vbox{\hrule\phantom{o}\hrule}\vrule}}

\def\nr{\par \noindent}

\def\beq{\begin{equation}}
\def\eeq{\end{equation}}

\def\ba{\begin{array}}
\def\ea{\end{array}}
\def\beann{\begin{eqnarray*}}
\def\eeann{\end{eqnarray*}}
\def\bea{\begin{eqnarray}}
\def\eea{\end{eqnarray}}
\def\Tr{\mathop{{\rm Trace}}}

\def\Argmin{\mathop{\rm Argmin}}

\def\Bbb#1{{\bf #1}}

\def\fnote#1{\footnote}
\def\blacksquare{\hbox{\vrule width 4pt height 4pt depth 0pt}}
\def\square{\hbox{\vrule\vbox{\hrule\phantom{o}\hrule}\vrule}}

\def\nr{\par \noindent}

\def\Tr{\mathop{{\rm Tr}}}

\def\Argmin{\mathop{\rm Argmin}}

\def\bar{\widehat}

\def\cB{{\cal B}}

\def\cN{{\cal N}}

\def\epsilonvi{{\epsilon_{\hbox{\scriptsize\rm VI}}}}

\def\Argmin{\mathop{\hbox{\rm  Argmin}}}

\def\Bbb#1{{\bf #1}}
\def\beq{\begin{equation}}

\def\blacksquare{\hbox{\vrule width 4pt height 4pt depth 0pt}}
\def\bR{{\mathbb{R}}}

\def\Card{{\hbox{\rm  Card}}}

\def\eeq{\end{equation}}

\def\fnote#1{\footnote}

\def\mypict3{\epsfxsize=220pt\epsfysize=80pt\epsfbox}

\def\nr{\par \noindent}

\def\Opt{{\hbox{\rm  Opt}}}

\def\Res{{\hbox{\rm  Res}}}

\def\square{\hbox{\vrule\vbox{\hrule\phantom{o}\hrule}\vrule}}

\def\Tr{{\hbox{\rm  Tr}}}

\def\cL{{{\cal L}}}

\def\bR{{\mathbf{R}}}
\def\Res{{\hbox{\rm Res}}}
\def\nuc{{\hbox{\scriptsize\rm nuc}}} 

%% file: MDLin2014-arXiv.bbl
\begin{thebibliography}{10}

\bibitem{Mosek}
E.~D. Andersen and K.~D. Andersen.
\newblock The {MOSEK} optimization tools manual.
\newblock
  \url{http://www.mosek.com/fileadmin/products/6_0/tools/doc/pdf/tools.pdf}.

\bibitem{Chambol05}
J.-F. Aujol and A.~Chambolle.
\newblock Dual norms and image decomposition models.
\newblock {\em International Journal of Computer Vision}, 63(1):85--104, 2005.

\bibitem{Teboulle09a}
A.~Beck and M.~Teboulle.
\newblock A fast iterative shrinkage-thresholding algorithm for linear inverse
  problems.
\newblock {\em SIAM Journal on Imaging Sciences}, 2(1):183--202, 2009.

\bibitem{R100}
S.~Becker, J.~Bobin, and E.~J. Cand{\`e}s.
\newblock Nesta: a fast and accurate first-order method for sparse recovery.
\newblock {\em SIAM Journal on Imaging Sciences}, 4(1):1--39, 2011.

\bibitem{Boyd10}
S.~Boyd, N.~Parikh, E.~Chu, B.~Peleato, and J.~Eckstein.
\newblock Distributed optimization and statistical learning via the alternating
  direction method of multipliers.
\newblock {\em Foundations and Trends® in Machine Learning}, 3(1):122--122,
  2010.

\bibitem{Buades05}
A.~Buades, B.~Coll, and J.-M. Morel.
\newblock A review of image denoising algorithms, with a new one.
\newblock {\em Multiscale Modeling \& Simulation}, 4(2):490--530, 2005.

\bibitem{Candes11}
E.~J. Cand\'{e}s, X.~Li, Y.~Ma, and J.~Wright.
\newblock Robust principal component analysis?
\newblock {\em Journal of the ACM (JACM)}, 58(3):11, 2011.

\bibitem{R101}
A.~Chambolle and T.~Pock.
\newblock A first-order primal-dual algorithm for convex problems with
  applications to imaging.
\newblock {\em Journal of Mathematical Imaging and Vision}, 40(1):120--145,
  2011.

\bibitem{CT93}
G.~Chen and M.~Teboulle.
\newblock Convergence analysis of a proximal-like minimization algorithm using
  bregman functions.
\newblock {\em SIAM Journal on Optimization}, 3(3):538--543, 1993.

\bibitem{Wo13}
W.~Deng, M.-J. Lai, Z.~Peng, and W.~Yin.
\newblock Parallel multi-block admm with o (1/k) convergence, 2013.
\newblock \url{http://www.optimization-online.org/DB_HTML/2014/03/4282.html}.

\bibitem{Goldfarb12}
D.~Goldfarb and S.~Ma.
\newblock Fast multiple-splitting algorithms for convex optimization.
\newblock {\em SIAM Journal on Optimization}, 22(2):533--556, 2012.

\bibitem{Goldfarb13}
D.~Goldfarb, S.~Ma, and K.~Scheinberg.
\newblock Fast alternating linearization methods for minimizing the sum of two
  convex functions.
\newblock {\em Mathematical Programming}, 141(1-2):349--382, 2013.

\bibitem{CVX}
M.~Grant, S.~Boyd, and Y.~Ye.
\newblock Cvx: Matlab software for disciplined convex programming, 2008.

\bibitem{FOM11}
A.~Juditsky and A.~Nemirovski.
\newblock First-order methods for nonsmooth largescale convex minimization: I
  general purpose methods; ii utilizing problems structure.
\newblock In S.~Sra, S.~Nowozin, and S.~Wright, editors, {\em Optimization for
  Machine Learning}, pages 121--183. The MIT Press, 2011.

\bibitem{LNN}
C.~Lemaréchal, A.~Nemirovskii, and Y.~Nesterov.
\newblock New variants of bundle methods.
\newblock {\em Mathematical Programming}, 69(1-3):111--147, 1995.

\bibitem{Monteiro13}
R.~D. Monteiro and B.~F. Svaiter.
\newblock Iteration-complexity of block-decomposition algorithms and the
  alternating direction method of multipliers.
\newblock {\em SIAM Journal on Optimization}, 23(1):475--507, 2013.

\bibitem{MP}
A.~Nemirovski.
\newblock Prox-method with rate of convergence o (1/t) for variational
  inequalities with lipschitz continuous monotone operators and smooth
  convex-concave saddle point problems.
\newblock {\em SIAM Journal on Optimization}, 15(1):229--251, 2004.

\bibitem{NOR}
A.~Nemirovski, S.~Onn, and U.~G. Rothblum.
\newblock Accuracy certificates for computational problems with convex
  structure.
\newblock {\em Mathematics of Operations Research}, 35(1):52--78, 2010.

\bibitem{NR02}
A.~Nemirovski and R.~Rubinstein.
\newblock An efficient stochastic approximation algorithm for stochastic saddle
  point problems.
\newblock In M.~Dror, P.~L'Ecuyer, and F.~Szidarovszky, editors, {\em Modeling
  Uncertainty and examination of stochastic theory, methods, and applications},
  pages 155--184. Kluwer Academic Publishers, 2002.

\bibitem{NesCompMin}
Y.~Nesterov.
\newblock Gradient methods for minimizing composite functions.
\newblock {\em Mathematical Programming}, 140(1):125--161, 2013.

\bibitem{Srebro12}
F.~Orabona, A.~Argyriou, and N.~Srebro.
\newblock Prisma: Proximal iterative smoothing algorithm.
\newblock {\em arXiv preprint arXiv:1206.2372}, 2012.

\bibitem{Lan2014}
Y.~Ouyang, Y.~Chen, G.~Lan, and E.~Pasiliao~Jr.
\newblock An accelerated linearized alternating direction method of
  multipliers, 2014.
\newblock \url{http://arxiv.org/abs/1401.6607}.

\bibitem{Qin12}
Z.~Qin and D.~Goldfarb.
\newblock Structured sparsity via alternating direction methods.
\newblock {\em The Journal of Machine Learning Research}, 13:1373--1406, 2012.

\bibitem{Katya11}
K.~Scheinberg, D.~Goldfarb, and X.~Bai.
\newblock Fast first-order methods for composite convex optimization with
  backtracking.
\newblock \url{http://www.optimization-online.org/DB_FILE/2011/04/3004.pdf},
  2011.
\newblock \url{http://www.optimization-online.org/DB_FILE/2011/04/3004.pdf}.

\bibitem{R102}
P.~Tseng.
\newblock Alternating projection-proximal methods for convex programming and
  variational inequalities.
\newblock {\em SIAM Journal on Optimization}, 7(4):951--965, 1997.

\bibitem{Tseng08}
P.~Tseng.
\newblock On accelerated proximal gradient methods for convex-concave
  optimization.
\newblock {\em submitted to SIAM Journal on Optimization}, 2008.

\bibitem{Goldfarb10}
Z.~Wen, D.~Goldfarb, and W.~Yin.
\newblock Alternating direction augmented lagrangian methods for semidefinite
  programming.
\newblock {\em Mathematical Programming Computation}, 2(3-4):203--230, 2010.

\end{thebibliography}
